\documentclass{gtart_h}

\def\ifplaintex{\expandafter\ifx\csname documentclass\endcsname\relax}

\ifplaintex 
\else
\fi

\def\gt{{\mathsurround=0pt\it $\cal G\mskip-2mu$eometry \&\ 
$\cal T\!\!$opology}}        

\def\gtm{{\mathsurround=0pt\it $\cal G\mskip-2mu$eometry \&\ 
$\cal T\!\!$opology $\cal M\mskip-1mu$onographs}}    

\def\gtp{{\mathsurround=0pt\it $\cal G\mskip-2mu$eometry \&\ 
$\cal T\!\!$opology $\cal P\!$ublications}}  

\def\recd{{\small Received:\qua\receiveddate\ifx\reviseddate\relax
\else\qquad Revised:\qua\reviseddate\fi\par}} 

\def\volumenumber#1{\def\thevolumenumber{#1}}
\def\volumeyear#1{\def\thevolumeyear{#1}}
\def\volumename#1{\def\thevolumename{#1}}
\def\papernumber#1{\def\thepapernumber{#1}}
\def\pagenumbers#1#2{\def\startpage{#1}\def\finishpage{#2}}
\def\published#1{\def\publishdate{#1}}
\def\received#1{\def\receiveddate{#1}}
\def\revised#1{\def\reviseddate{#1}}
\def\accepted#1{\def\accepteddate{#1}}

\def\coverauthors#1{\def\thecoverauthors{#1}}
\def\asciiauthors#1{\def\theasciiauthors{#1}}

\def\coverauthors#1{\def\thecoverauthors{#1}}
\long\def\asciiabstract#1{\long\def\theasciiabstract{#1}}

\let\\\par
\let\thevolumenumber\relax\let\thepapernumber\relax
\let\thevolumeyear\relax\let\startpage\relax
\let\finishpage\relax\let\publishdate\relax\let\receiveddate\relax
\let\reviseddate\relax\let\accepteddate\relax\let\theasciititle\relax
\let\theasciiauthors\relax
\let\theasciiabstract\relax
\let\thecoverauthors\relax
\let\thecoverauthors\relax\let\theerratum\relax\let\theasciiemail\relax
\let\theshortauthors\relax\let\theshorttitle\relax

\def\startpage{1}\def\finishpage{15}\def\thepapernumber{77}

\volumenumber{2}
\volumename{Proceedings of the Kirbyfest}
\volumeyear{1999}

\long\def\maketitlep{   

\count0=\startpage

\gtm\nl        
{\small Volume \thevolumenumber: \thevolumename\nl 
\ifx\theerratum\relax\else Erratum \erratumnumber\nl\fi
Pages \startpage--\finishpage\nl}

\vglue 0.1truein   

{\parskip=0pt\leftskip 0pt plus 1fil\def\\{\par\smallskip}{\ifplaintex\large
\else\Large\fi\bf\thetitle}\par\medskip}   
\vglue 0.05truein 

{\parskip=0pt\leftskip 0pt plus 1fil\def\\{\par}{\sc\theauthors}
\par\medskip}%
 
\vglue 0.03truein 

{\small\leftskip 25pt\rightskip 25pt{\bf Abstract}\stdspace\theabstract

{\bf AMS Classification}\stdspace\theprimaryclass
\ifx\thesecondaryclass\relax\else; \thesecondaryclass\fi\par
{\bf Keywords}\stdspace \thekeywords\par}\vglue 7pt

}   

\font\phead=cmsl9 scaled 950
\font=cmsl9 scaled 1050
\font\pnum=cmbx10 scaled 913
\font\lnum=cmbx10 
\font\pfoot=cmsl9 scaled 950
\font=cmsl9 scaled 1050
\ifplaintex
\headline{\vbox to 0pt{\vskip -4.5mm\line{\small\phead\ifnum
\count0=\startpage ISSN 1464-8997 (on line)
1464-8989 (printed) \hfill {\pnum\folio}\else\ifodd\count0\def\\{ }%
\ifx\theshorttitle\relax\thetitle\else\theshorttitle\fi\hfill{\pnum\folio}
\else\def\\{ and }{\pnum\folio}\hfill\ifx\theshortauthors\relax\theauthors
\else\theshortauthors\fi\fi\fi}\vss}}
\footline{\vbox to 0pt{\vglue 0mm\line{\small\pfoot\ifnum\count0=\startpage
Published \publishdate:\qua\copyright\ \gtp\hfill\else
\gtm, Volume \thevolumenumber\ (\thevolumeyear)\hfill\fi}\vss
}}
\else
\makeatletter
\def\@oddhead{{\small\lhead\ifnum\count0=\startpage ISSN 1464-8997 (on line)
1464-8989 (printed) \hfill {\lnum\number\count0}\else\ifodd\count0
\def\\{ }\ifx\theshorttitle\relax \thetitle \else\theshorttitle\fi\hfill
{\lnum\number\count0}\else\def\\{ and }{\lnum\number\count0}
\hfill\ifx\theshortauthors\relax 
\theauthors\else\theshortauthors\fi\fi\fi}}\def\@evenhead{@oddhead}
\def\@oddfoot{\small\lfoot\ifnum\count0=\startpage Published \publishdate:\qua\copyright\ \gtp\hfill\else
\gtm, Volume \thevolumenumber\ (\thevolumeyear)\hfill\fi}
\def\@evenfoot{@oddfoot}
\makeatother
\fi

\let\maketitlepage\maketitlep
\let\makeshorttitle\maketitlepage
\let\maketitle\maketitlepage

\newwrite\gtoutfile
\long\gdef\makeheadfile{  
{\def\\{, }\def\s{ }
\immediate\openout\gtoutfile head.xxx
\immediate\write\gtoutfile{Proxy-for: \ifx\theasciiauthors\relax
\theauthors\else\theasciiauthors\fi\s<\ifx\theasciiemail\relax\theemail\else\theasciiemail\fi>}
\immediate\write\gtoutfile{\noexpand\\}
\immediate\write\gtoutfile{Authors: \ifx\theasciiauthors\relax
\theauthors\else\theasciiauthors\fi}
{\def\\{ }\immediate\write\gtoutfile{Title: \ifx\theasciititle\relax
\thetitle\else\theasciititle\fi}}
\immediate\write\gtoutfile{Subj-class: GT or SG, GR etc}
\immediate\write\gtoutfile{MSC-class: \theprimaryclass\ifx\thesecondaryclass\relax\else, \thesecondaryclass\fi}
\immediate\write\gtoutfile{Journal-ref: Geom. Topol. Monogr. \thevolumenumber\s
(\thevolumeyear) \startpage-\finishpage}
\immediate\write\gtoutfile{Comments: Published by Geometry and Topology Monographs at}
\immediate\write\gtoutfile{\s\s\s  http://www.maths.warwick.ac.uk/gt/GTMon\thevolumenumber/paper\thepapernumber.abs.html}
\immediate\write\gtoutfile{\noexpand\\}
\immediate\write\gtoutfile{}
\ifx\theasciiabstract\relax
\immediate\write\gtoutfile{\theabstract}\else
\immediate\write\gtoutfile{\theasciiabstract}\fi
\immediate\write\gtoutfile{}
\immediate\write\gtoutfile{\noexpand\\}
\immediate\write\gtoutfile{}
\immediate\closeout\gtoutfile}}  

\def\maketitlepage{\maketitlep\makeheadfile}
\let\makeshorttitle\maketitlepage
\let\maketitle\maketitlepage

\volumenumber{7}
\volumename{Proceedings of the Casson Fest}
\volumeyear{2004}
\papernumber{2}
\pagenumbers{27}{68}
\received{10 November 2003}
\revised{10 March 2004}
\accepted{10 March 2004}
\published{17 September 2004}

\usepackage{amsmath,amssymb}
\usepackage{graphicx}

\newtheorem{thm}{Theorem}[section]
\newtheorem{lem}[thm]{Lemma}
\newtheorem{cor}[thm]{Corollary}

\newtheorem{case}{Case}
\newtheorem{claim*}{Claim}

\theoremstyle{remark}
\newtheorem*{rem}{Remark}

\newcommand{\ben}{\begin{enumerate}}
\newcommand{\een}{\end{enumerate}}

\newcommand{\wh}{\widehat}
\newcommand{\wt}{\widetilde}
\newcommand{\mc}{\mathcal}
\newcommand{\intr}{\text{int}\,}

\newcommand{\set}[1]{\left\{#1\right\}}

\newcommand{\sk}{\mc{SK}}

\newcommand{\Figw}[4]{
\begin{center}
\includegraphics[width=#1]{#2}
\end{center}
\caption{ #3 \label{#4} } }

\begin{document}

\title{Seifert Klein bottles for knots with\\common boundary slopes}

\author{Luis G Valdez-S\'anchez}
\coverauthors{Luis G Valdez-S\noexpand\'anchez}
\asciiauthors{Luis G Valdez-Sanchez}

\address{Department of Mathematical Sciences,
University of Texas at El Paso\\
El Paso, TX 79968, USA}

\email{valdez@math.utep.edu}

\begin{abstract}   
We consider the question of how many essential Seifert Klein
bottles with common boundary slope a knot in $S^3$ can bound, up
to ambient isotopy. We prove that any hyperbolic knot in $S^3$
bounds at most six Seifert Klein bottles with a given boundary
slope. The Seifert Klein bottles in a minimal projection of
hyperbolic pretzel knots of length 3 are shown to be unique and
$\pi_1$--injective, with surgery along their boundary slope
producing irreducible toroidal manifolds. The cable knots which
bound essential Seifert Klein bottles are classified; their
Seifert Klein bottles are shown to be non-$\pi_1$--injective, and
unique in the case of torus knots. For satellite knots we show
that, in general, there is no upper bound for the number of
distinct Seifert Klein bottles a knot can bound.
\end{abstract}

\asciiabstract{%
We consider the question of how many essential Seifert Klein
bottles with common boundary slope a knot in S^3 can bound, up
to ambient isotopy. We prove that any hyperbolic knot in S^3
bounds at most six Seifert Klein bottles with a given boundary
slope. The Seifert Klein bottles in a minimal projection of
hyperbolic pretzel knots of length 3 are shown to be unique and
pi_1-injective, with surgery along their boundary slope
producing irreducible toroidal manifolds. The cable knots which
bound essential Seifert Klein bottles are classified; their
Seifert Klein bottles are shown to be non-pi_1-injective, and
unique in the case of torus knots. For satellite knots we show
that, in general, there is no upper bound for the number of
distinct Seifert Klein bottles a knot can bound.}


\primaryclass{57M25}
\secondaryclass{57N10}
\keywords{Seifert Klein bottles, knot complements, boundary slope}

\maketitle
\cl{\small\it Dedicated to Andrew J Casson on the
occasion of his 60th birthday}

\section{Introduction}\label{intro}

For any knot in $S^3$ all orientable Seifert surfaces spanned by
the knot have the same boundary slope. The smallest genus of such
a surface is called the {\it genus} of the knot, and such a
minimal surface is always essential in the knot exterior.
Moreover, by a result of Schubert--Soltsien \cite{schubert1}, any
simple knot admits finitely many distinct minimal genus Seifert
surfaces, up to ambient isotopy, while for satellite knots
infinitely many isotopy classes may exist (cf
\cite{eisner1,kakimizu1}).

The smallest genus of the nonorientable Seifert surfaces spanned
by a knot is the {\it crosscap number} of the knot (cf
\cite{clark}). Unlike their orientable counterparts, nonorientable
minimal Seifert surfaces need not have a unique boundary slope. In
fact, by \cite{tera1}, any knot $K\subset S^3$ has at most two
boundary slopes $r_1,r_2$ corresponding to {\it essential}
(incompressible and boundary incompressible, in the geometric
sense) Seifert Klein bottles, and if so then $\Delta(r_1,r_2)=4$
or $8$, the latter distance occurring only when $K$ is the
figure--8 knot. The knots for which two such slopes exist were
classified in \cite{valdez6} and, with the exception of the
$(2,1,1)$ and $(-2,3,7)$ pretzel knots (the figure--8 knot and the
Fintushel--Stern knot, respectively), all are certain satellites of
2--cable knots. Also, a minimal crosscap number Seifert surface for
a knot need not even be essential in the knot exterior (cf
\cite{tera1}).

In this paper we study the uniqueness or non-uniqueness, up to
ambient isotopy, of essential Seifert Klein bottles for a knot
with a fixed boundary slope; we will regard any two such surfaces
as {\it equivalent} iff they are ambient isotopic. Our main result
states that any crosscap number two hyperbolic knot admits at most
6 nonequivalent Seifert Klein bottles with a given slope; before
stating our results in full we will need some definitions.

We work in the PL category; all 3--manifolds are assumed to be
compact and orientable. We refer to ambient isotopies simply as
isotopies. Let $M^3$ be a 3--manifold with boundary. The pair
$(M^3,\partial M^3)$ is {\it irreducible} if $M^3$ is irreducible
and $\partial M^3$ is incompressible in $M^3$. Any embedded circle
in a once punctured Klein bottle is either a {\it meridian}
(orientation preserving and nonseparating), a {\it longitude}
(orientation preserving and separating), or a {\it center}
(orientation reversing); any two meridians are isotopic within the
surface, but there are infinitely many isotopy classes of
longitudes and centers (cf Lemma~\ref{circles}). For a knot $K$
in $S^3$ with exterior $X_K=S^3\setminus\intr N(K)$ and a
nontrivial slope $r$ in $\partial X_K$, $K(r)=X_K(r)$ denotes the
manifold obtained by Dehn--filling $X_K$ along $r$, that is, the
result of {\it surgering $K$ along $r$}. We denote by $\sk(K,r)$
the collection of equivalence classes of essential Seifert Klein
bottles in $X_K$ with boundary slope $r$; as pointed out above,
$\sk(K,r)$ is nonempty for at most two distinct integral slopes
$r_1,r_2$. If $|\sk(K,r)|\geq 2$, we will say that the collection
$\sk(K,r)$ is {\it meridional} if any two distinct elements can be
isotoped so as to intersect transversely in a common meridian, and
that it is {\it central} if there is a link $c_1\cup c_2$ in $X_K$
such that any two distinct elements of $\sk(K,r)$ can be isotoped
so as to intersect transversely in $c_1\cup c_2$, and $c_1,c_2$
are disjoint centers in each element. For $P\in\sk(K,r)$, $N(P)$
denotes a small regular neighborhood of $P$, and
$H(P)=X_K\setminus\intr N(P)$ denotes the {\it exterior} of $P$ in
$X_K$. We say $P$ is {\it unknotted} if $H(P)$ is a handlebody,
and {\it knotted} otherwise; if the pair $(H(P),\partial H(P))$ is
irreducible, we say $P$ is {\it strongly knotted}. If
$\mu,\lambda$ is a standard meridian--longitude pair for  a knot
$L$, and  $K$ is a circle embedded in $\partial X_L$ representing
$p\mu+q\lambda$ for some relatively prime integers $p,q$ with
$|q|\geq 2$, we say $K$ is a {\it $(p,q)$ cable of $L$}; we also
call $K$ a {\it $q$--cable knot}, or simply a {\it cable knot}. In
particular, the torus knot $T(p,q)$ is the $(p,q)$--cable of the
trivial knot. If $X$ is a finite set or a topological space, $|X|$
denotes its number of elements or of connected components.

\begin{thm}\label{t1}
Let $K$ be a hyperbolic knot in $S^3$ and $r$  a slope in
$\partial X_K$ for which $\sk(K,r)$ is nonempty. Then any element
of $\sk(K,r)$ is either unknotted or strongly knotted, and

\ben
\item[\rm(a)] if $|\sk(K,r)|\geq 2$ then $\sk(K,r)$ is
either central or meridional; in the first case the link $c_1\cup
c_2$ is unique up to isotopy in $X_K$ and $|\sk(K,r)|\leq 3$, in
the latter case $|\sk(K,r)|\leq 6$;

\item[\rm(b)] if $\sk(K,r)$ is central then each of its elements is
strongly knotted; if some element of $\sk(K,r)$ is unknotted then
$|\sk(K,r)|\leq 2$, and $\sk(K,r)$ is meridional if
$|\sk(K,r)|=2$;

\item[\rm(c)]
if some element of $\sk(K,r)$ is not $\pi_1$--injective then it is
unknotted, $K$ has tunnel number one, and $K(r)$ is a Seifert
fibered space over $S^2$ with at most 3 singular fibers of indices
$2,2,n$ and finite fundamental group;

\item[\rm(d)] if some element of $\sk(K,r)$ is $\pi_1$--injective and
unknotted then $K(r)$ is irreducible and toroidal.
\een
\end{thm}

\begin{cor}
Let $K\subset S^3$ be a hyperbolic knot and $P$ an unknotted
element of $\sk(K,r)$. Then $\pi_1(K(r))$ is finite iff $P$ is not
$\pi_1$--injective. In particular, if $r=0$, then $P$ is
$\pi_1$--injective.\hfill\qed
\end{cor}

In the case when $\sk(K,r)$ is meridional and contains an
unknotted element we give examples in Section \ref{submer} realizing the
bound $|\sk(K,r)|=2$ for $K$ a hyperbolic knot. Such examples are
obtained from direct variations on the knots constructed by Lyon
in \cite{lyon1}; M. Teragaito (personal communication) has constructed
more examples along similar lines. It is not known if the other
bounds given in Theorem~\ref{t1} are optimal, but see the remark
after Lemma~\ref{four} for a discussion on possible ways of
realizing the bound $|\sk(K,r)|=4$, and Section \ref{subcen} for a
construction of possible examples of central families $\sk(K,r)$
with $|\sk(K,r)|=2$. On the other hand, hyperbolic knots which
span a unique Seifert Klein bottle per slope are not hard to find:
if we call {\it algorithmic} any black or white surface obtained
from some regular projection of a knot, then the algorithmic
Seifert Klein bottles in minimal projections of hyperbolic pretzel
knots provide the simplest examples.

\begin{thm}\label{t2}
Let $K\subset S^3$ be a hyperbolic pretzel knot, and let $P$ be
any algorithmic Seifert Klein bottle in a minimal projection of
$K$. Then $P$ is unknotted, $\pi_1$--injective, and unique up to
equivalence; moreover, $K(\partial P)$ is irreducible and
toroidal.
\end{thm}

We remark that any crosscap number two 2--bridge knot is a
hyperbolic pretzel knot. Also, as mentioned before, the only
hyperbolic knots which bound Seifert Klein bottles of distinct
boundary slope are the $(2,1,1)$ and $(-2,3,7)$ pretzel knots. The
standard projection of the $(2,1,1)$ pretzel knot simultaneously
realizes two algorithmic Seifert Klein bottles of distinct slopes
(in fact, this is the only nontrivial knot in $S^3$ with such property), so
they are handled by Theorem~\ref{t2}. The $(-2,3,7)$ pretzel knot
has both an algorithmic and a non algorithmic Seifert Klein
bottle; it can be proved that the non algorithmic surface also
satisfies the conclusion of Theorem~\ref{t2}, but we omit the
details.

In contrast, no universal bound for $|\sk(K,r)|$ exists for
satellite knots:

\begin{thm}\label{t3}
For any positive integer $N$, there are satellite knots $K\subset
S^3$ with $|\sk(K,0)|\geq N$.
\end{thm}

Among non hyperbolic knots, the families of cable or composite
knots are of particular interest; we classify the crosscap number
two cable knots, and find information about the Seifert Klein
bottles bounded by composite knots.

\begin{thm}\label{t4}
Let $K$ be a knot in $S^3$ whose exterior contains an essential
annulus and an essential Seifert Klein bottle. Then either
\ben
\item[\rm(a)] $K$ is a $(2(2m+1)n\pm 1,4n)$--cable for some integers
$m,n,  \ n\neq 0$,

\item[\rm(b)] $K$ is a $(2(2m+1)(2n+1)\pm 2,2n+1)$--cable of a
$(2m+1,2)$--cable knot for some integers $m,n$,

\item[\rm(c)] $K$ is a $(6(2n+1)\pm 1,3)$--cable of some $(2n+1,2)$--cable knot
for some integer $n$, or

\item[\rm(d)] $K$ is a connected sum of two 2--cable knots.
\een
Any Seifert Klein bottle bounded by any composite knot is
$\pi_1$--injective, while the opposite holds for any cable knot; in
case (a), such a surface is unique up to equivalence.
\end{thm}

The next result follows from the proof of Theorem~\ref{t4}:

\begin{cor}\label{torus}
The crosscap number two torus knots are $T(\pm 5,3)$, $T(\pm
7,3)$, and $T(2(2m+1)n\pm 1,4n)$ for some $m,n,  \ n\neq 0$; each
bounds a unique Seifert Klein bottle, which is unknotted and not
$\pi_1$--injective.
\end{cor}

The crosscap numbers of torus knots have been determined in
\cite{tera4}; our classification of these knots follows
from Theorem~\ref{t4}, whose proof gives detailed topological
information about the construction of Seifert Klein bottles for
cable knots in general. The classification of crosscap number two
composite knots also follows from \cite{tera5} (where
many-punctured Klein bottles are considered); these knots serve as
examples of satellite knots any of whose Seifert Klein bottles is
$\pi_1$--injective. The Seifert Klein bottles  for the knots in
Theorem~\ref{t4}(a),(b) are disjoint from the cabling annulus, and
in (b) one expects the number of Seifert Klein bottles bounded by
the knot to depend on the nature of its companions, as in
\cite[Corollary D]{kakimizu2}; in (c),(d) the Seifert Klein
bottles intersect the cabling or splitting annulus, so in (d) one
would expect there to be infinitely nonequivalent Seifert Klein
bottles bounded by such knots, as in \cite{eisner1}.

The paper is organized as follows. Section~\ref{prelim} collects a
few more definitions and some general properties of Seifert Klein
bottles. In Section~\ref{satknots} we look at a certain family of
crosscap number two satellite knots which bound Seifert Klein
bottles with zero boundary slope, and use it to prove
Theorem~\ref{t3}. In Section~\ref{cable-comp} we first identify
the minimal intersection between an essential annulus and an
essential Seifert Klein bottle in a knot exterior, which is the
starting point of the proof of Theorem~\ref{t4} and its corollary.
Section~\ref{handlebodies} contains some results on non boundary
parallel separating annuli and pairs of pants contained in
3--manifolds with boundary, from both algebraic and geometric
points of view. These results have direct applications to the case
of unknotted Seifert Klein bottles, but we will see in
Lemma~\ref{hbdy-irred} that if $K$ is any hyperbolic knot and
$P,Q\in\sk(K,r)$ are any two distinct elements, then $Q$ splits
the exterior $H(P)$ of $P$ into two pieces, at least one of which
is a genus two handlebody; this observation eventually leads to
the proof of Theorem~\ref{t1} in Section~\ref{hypknots}. After
these developments, a proof of Theorem~\ref{t2} is given within a
mostly algebraic setting in Section~\ref{pretzels}.

\sh{Acknowledgements}

I want to thank  Francisco Gonz\'alez-Acu\~na, Enrique Ram\'{\i}rez-Losada
and V\'{\i}ctor N\'u\~nez for many useful conversations, and Masakazu
Teragaito and the Hiroshima University Department of Mathematics for their
hospitality while visiting them in the Spring of 2003, when part of this
work was completed. Finally, I would like to thank the referee for his
many useful suggestions and corrections to the original text.

\section{Preliminaries}\label{prelim}

In this section we set some more notation we will use in the
sequel, and establish some general properties of essential Seifert
Klein bottles. Let $M^3$ be a 3--manifold with boundary. For any
surface $F$ properly embedded in $M^3$ and  $c$ the union of some
components of $\partial F$, $\wh{F}$ denotes the surface in
$M^3(c)=M^3\cup\set{\text{2--handles along }c}$ obtained by capping
off the circles of $\partial F$ isotopic to $c$ in $\partial M^3$
suitably with disjoint disks in $M^3(c)$. If $G$ is a second
surface properly embedded in $M^3$ which intersects $F$
transversely with $\partial F\cap\partial G=\emptyset$, let
$N(F\cap G)$ be a small regular neighborhood of $F\cap G$ in
$M^3$, and let $\mc{A}$ be the collection of annuli obtained as
the closures of the components of $\partial N(F\cap
G)\setminus(F\cup G)$. We use the notation $F\asymp G$ to
represent the surface obtained by capping off the boundary
components in $\intr M^3$ of $F\cup G\setminus\intr N(F\cap G)$
with suitable annuli from $\mc{A}$. As usual,
$\Delta(\alpha,\beta)$ denotes the minimal geometric intersection
number between circles of slopes $\alpha,\beta$ embedded in a
torus.

Now let $P$ be a Seifert Klein bottle for a knot $K\subset S^3$,
and let $N(P)$ be a small regular neighborhood of $P$ in $X_K$;
$N(P)$ is an $I$--bundle over $P$, topologically a genus two
handlebody. Let $H(P)=X_K\setminus\intr N(P)$ be the exterior of
$P$ in $X_K$. Let $A_K,A_K'$ denote the annuli $N(P)\cap\partial
X_K,H(P)\cap\partial X_K$, respectively, so that $A_K\cup
A_K'=\partial X_K$ and $A_K\cap A_K'=\partial A_K=\partial A_K'$.
Then $\partial P$ is a core of $A_K$, and we denote the core of
$A_K'$ by $K'$. Finally, let $T_P$ denote the frontier of $N(P)$
in $X_K$. $T_P$ is a twice-punctured torus such that $N(P)\cap
H(P)=T_P$; since $N(P)$ is an $I$--bundle over $P$, $P$ is
$\pi_1$--injective in $N(P)$ and $T_P$ is incompressible in $N(P)$.

For any meridian circle $m$ of $P$, there is an annulus $A(m)$
properly embedded in $N(P)$ with $P\cap A(m)=m$ and $\partial
A(m)\subset T_P$; we call the circles $\partial A(m)=m_1\cup m_2$
the {\it lifts of $m$ (to $T_P$)}. Similarly, for any center
circle $c$ of $P$, there is a Moebius band $B(c)$ properly
embedded in $N(P)$ with $P\cap B(c)=c$ and $\partial B(c)\subset
T_P$; we call $l=\partial B(c)$ the {\it lift of $c$ (to $T_P$)}.
For a pair of disjoint centers in $P$, similar disjoint Moebius
bands can be found in $N(P)$. Since the meridian of $P$ is unique
up to isotopy, the lifts of a meridian of $P$ are also unique up
to isotopy in $T_P$; the lift of a center circle of $P$ depends
only on the isotopy class of the center circle in $P$.

We denote the linking form in $S^3$ by $\ell k(\cdot,\cdot)$.

\begin{lem}\label{slope}
Let $P$ be a Seifert Klein bottle for a knot $K\subset S^3$. If
$m$ is the meridian circle of $P$, then the boundary slope of $P$
is $\pm \: 2 \: \ell k(K,m)$.
\end{lem}

\begin{proof}
For $m$ a meridian circle of $P$ let $P'$ be the pair of pants
$P\setminus\intr N(m)$ and $A'$ the annulus $P\cap N(m)$, where
$N$ is a small regular neighborhood of $m$ in $X_K$; thus
$P=P'\cup A'$, and $\partial P'=\partial P\cup
\partial A'$. Fixing an orientation of $P'$ induces an
orientation on $\partial P'$ such that the circles $\partial A'$
become coherently oriented in $A'$. In this way an orientation on
$m$ is induced, coherent with that of $\partial A'$, such that
$\ell k(K,m_1)=\ell k(K,m_2)=\ell k(K,m)$. As the slope of
$\partial P$ is integral, hence equal to $\pm\ell k(K,\partial
P)$, and $\ell k(K,m_1\cup m_2\cup\partial P)=0$, the lemma
follows.
\end{proof}

\begin{lem}\label{pinc}
If $P\in\sk(K,r)$ then $\wh{P}$ is incompressible in $K(r)$.
\end{lem}

\begin{proof}
If $\wh{P}$ compresses in $K(r)$ along a circle $\gamma$ then
$\gamma$ must be orientation-preserving in $\wh{P}$. Thus,
surgering $\wh{P}$ along a compression disk $D$ with $\partial
D=\gamma$ produces either a nonseparating 2--sphere (if $\gamma$
is a meridian) or two disjoint projective planes (if $\gamma$ is a
longitude) in $K(r)$, neither of which is possible as $K$
satisfies Property R \cite{gabai} and $K(r)$ has cyclic integral
homology.
\end{proof}

\begin{lem}\label{lifts}
Let $m$ be a meridian circle and $c_1,c_2$ be two disjoint center
circles of $P$; let $m_1,m_2$ and $l_1,l_2$ be the lifts of $m$
and $c_1,c_2$, respectively. Then,
\ben
\item[\rm(a)] neither circle $K',l_1,l_2$  bounds a surface in $H(P)$,

\item[\rm(b)] neither pair $m_1,m_2$ nor
$l_1,l_2$ cobound a surface in $H(P)$,

\item[\rm(c)] none of the circles $m_i,l_i$ cobounds an annulus in $H(P)$
with $K'$, and

\item[\rm(d)] if $A$ is an annulus in $H(P)$ with $\partial A=\partial
A_{K}'$, then $A$ is not parallel in $X_K$ into $A_K$.
\een
\end{lem}

\begin{proof}
Let $B_i$ be a Moebius band in $N(P)$ bounded by $l_i$. If $K'$ or
$l_i$ bounds a surface $F$ in $H(P)$ then out of the surfaces
$P,B_i,F$ it is possible to construct a nonorientable closed
surface in $S^3$, which is impossible; thus (a) holds.

Consider the circles $K',\alpha_1,\alpha_2$ in $\partial H(P)$,
where $\alpha_1,\alpha_2=m_1,m_2$ or $l_1,l_2$; such circles are
mutually disjoint. Let $P_1$ be the closure of some component of
$\partial H(P)\setminus(K'\cup\alpha_1\cup\alpha_2)$; $P_1$ is a
pair of pants. If $\alpha_1,\alpha_2$ cobound a surface $F$
properly embedded in $H(P)$ then $F\cup_{\alpha_1\cup\alpha_2}P_1$
is a surface in $H(P)$ bounded by $K'$, which can not be the case
by (a). Hence $\alpha_1,\alpha_2$ do not cobound a surface in
$H(P)$ and so (b) holds.

If any circle $m_1,m_2,l_1,l_2$ cobounds an annulus in $H(P)$ with
$K'$ then $\wh{P}$ compresses in $K(\partial P)$, which is not the
case by Lemma~\ref{pinc}; thus (c) holds.

Finally, if $A\subset H(P)$ is parallel to $A_K$ then the region
$V$ cobounded by $A\cup A_K$ is a solid torus with $P\subset V$
and $\partial P$ a core of $A_K$, hence $V(\partial P)=S^3$
contains the closed Klein bottle $\wh{P}$, which is impossible.
This proves (d).
\end{proof}

The next result follows from \cite[Theorem 1.3]{gordonlu5} and
\cite[Lemma 4.2]{valdez6}:

\begin{lem}\label{kleiness}
Let $K$ be a nontrivial knot in $S^3$. If $P$ is a once-punctured
Klein bottle properly embedded in $X_K$ then $P$ is essential iff
$K$ is not a 2--cable knot, and in such case $P$ has integral
boundary slope.\hfill\qed
\end{lem}

\section{The size of $\sk(K,0)$}\label{satknots}

In this section we consider a special family of crosscap number
two satellite knots $\{\mc{K}_n\}$, which generalizes the example
of W.R.\ Alford in \cite{alford1}; the knots in this family are
constructed as follows. For each $i\geq 1$ let $K_i$ be a
nontrivial {\it prime} knot in $S^3$. Figure~\ref{fig01} shows a
pair of pants $F_n$, $n\geq 2$, with only one shaded side,
constructed with a band $B_n$ whose core `follows the pattern' of
$K_1\# K_2\#\cdots\# K_n$, and which is attached to two disjoint,
unknotted, and untwisted annuli $A',A''$; we define the knot
$\mc{K}_n$ to be the boundary component of $F_n$ indicated in the
same figure.

\begin{figure}
\Figw{3.5in}{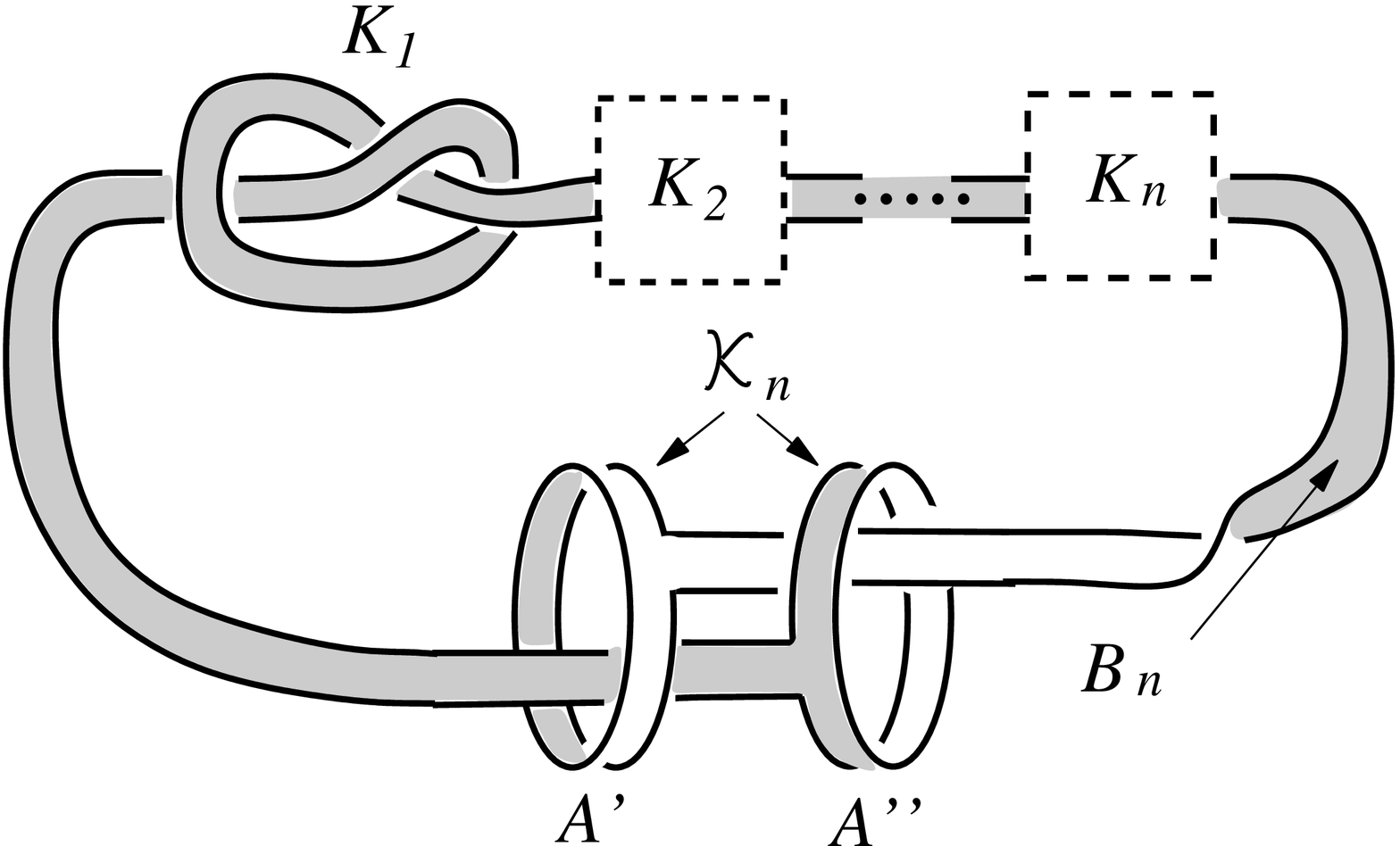}{The pair of pants $F_n=A'\cup B_n\cup
A''$ and the knot $\mc{K}_n\subset\partial F_n$}{fig01}
\end{figure}

For each $1\leq s< n$, let $P_s$ be the Seifert Klein bottle
bounded by $\mc{K}_n$ constructed by attaching an annulus $A_s$ to
the two boundary circles of $F_n$ other than $\mc{K}_n$, which
{\it swallows} the factors $K_1,\dots,K_s$ and {\it follows} the
factors $K_{s+1},\dots,K_n$, as indicated in Figure~\ref{fig02}.

Notice that the core of $A_s$, which has linking number zero with
$\mc{K}_n$, is a meridian circle of $P_s$, so the boundary slope
of $P_s$ is zero by Lemma~\ref{slope}. Our goal is to show that
$P_r$ and $P_s$ are not equivalent for $r\neq s$, so that
$|\sk(\mc{K}_n,0)|\geq n-1$, which will prove Theorem~\ref{t3}. In
fact, we will prove the stronger statement that $P_r$ and $P_s$
are not equivalent for $r\neq s$ even under homeomorphisms of
$S^3$.

\begin{figure}
\Figw{3.6in}{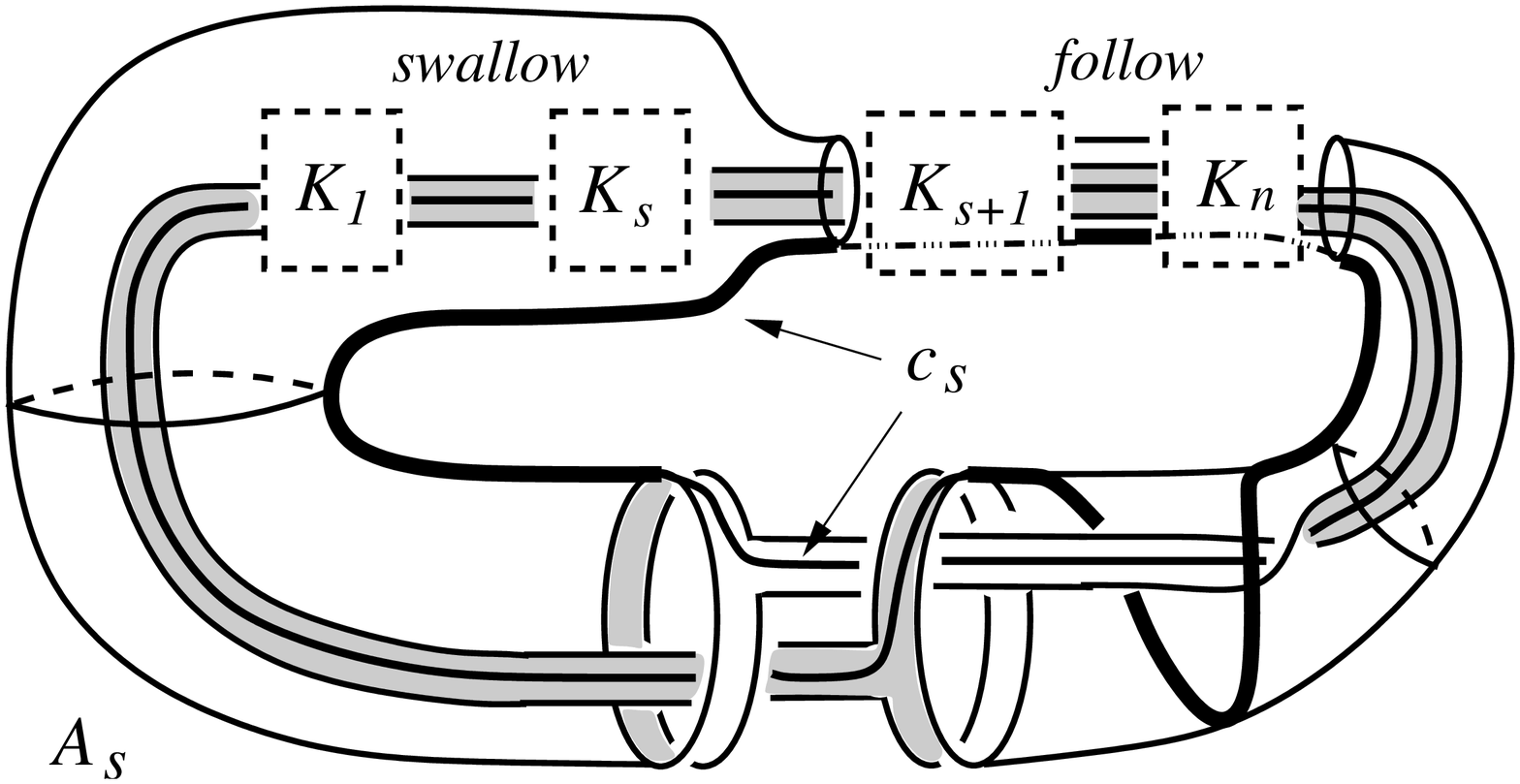}{The Seifert Klein bottle $P_s=F_n\cup A_s$
and center $c_s\subset P_s$}{fig02}
\end{figure}

The following elementary result on intersection properties between
essential circles in a once-punctured Klein bottle will be useful
in determining all centers of any of the above Seifert Klein
bottles of $\mc{K}_n$; we include its proof for the convenience of
the reader.

\begin{lem}\label{circles}
Let $P=A\cup B$ be a once-punctured Klein bottle, where $A$ is an
annulus and $B$ is a rectangle with $A\cap B=\partial
A\cap\partial B$ consisting of two opposite edges of $\partial B$,
one in each component of $\partial A$. Let $m$ be a core of $A$
(ie a meridian of $P$) and $b$ a core arc in $B$ parallel to
$A\cap B$. If $\omega$ is any nontrivial circle embedded in $P$
and not parallel to $\partial P$ which has been isotoped so as to
intersect $m\cup b$ transversely and minimally, then $\omega$ is a
meridian, center, or longitude circle of $P$ iff $(|\omega\cap
m|,|\omega\cap b|)=(0,0),\ (1,1)$, or $(2,2)$,
respectively.
\end{lem}

\begin{proof}
Let $I_1,I_2$ denote the components of $A\cap B=\partial
A\cap\partial B$. After isotoping  $\omega$  so as to intersect
$m\cup b$ transversely and minimally, either $\omega$ lies in
$\intr A$ and is parallel to $m$, or $\omega\cap B$  consists of
disjoint spanning arcs of $B$ with endpoints in $I_1\cup I_2$,
while $\omega\cap A$ consists of disjoint arcs which may split
into at most 4 parallelism classes, denoted
$\alpha,\beta,\gamma,\delta$; the situation is represented in
Figure~\ref{fig04}.

Suppose we are in the latter case. As $\omega$ is connected and
necessarily $|\alpha|=|\delta|$, if $|\alpha|>0$ then $\omega\cap
A$ must consist of one arc of type $\alpha$ and one arc of type
$\delta$; but then $\omega$ is parallel to $\partial P$, which is
not the case. Thus we must have $|\alpha|=|\delta|=0$, in which
case $|\beta|+|\gamma|=|\omega\cap b|=n$ for some integer $n\geq
1$. Notice that $\omega$ is a center if $n=1$ and a longitude if
$n=2$; in the first case  $\set{|\beta|,|\gamma|}=\set{0,1}$,
while in the latter $\set{|\beta|,|\gamma|}=\set{0,2}$. These are
the only possible options for $\omega$ and $n$ whenever
$|\beta|=0$ or $|\gamma|=0$.

Assume that $|\beta|,|\gamma|\geq 1$, so $n\geq 3$, and label the
endpoints of the arcs $\omega\cap A$ and $\omega\cap B$ in
$I_1,I_2$ consecutively with $1,2,\dots,n$, as in
Figure~\ref{fig04}. We assume, as we may, that
$|\beta|\leq|\gamma|$. We start traversing $\omega$ from the point
labelled $1$ in $I_1\subset\partial B$ in the direction of
$I_2\subset\partial B$, within $B$, then reach the endpoint of an
arc component of $\omega\cap A$ in $I_2\subset\partial A$, then
continue within $A$ to an endpoint in $I_1\subset\partial A$, and
so on until traversing all of $\omega$. The arc components of
$\omega\cap B$ and $\omega\cap A$ traversed in this way give rise
to permutations $\sigma,\tau$ of $1,2,\dots,n$, respectively,
given by $\sigma(x)=n-x+1$ for $1\leq x\leq n$, and
$\tau(x)=n+x-|\beta|$ for $1\leq x\leq|\beta|$,
$\tau(x)=x-|\beta|$ for $|\beta|<x\leq n$.

Clearly, the number of components of $\omega$ equals the number of
orbits of the permutation $\tau\circ\sigma$. But, since
$|\beta|\leq|\gamma|$, the orbit of $\tau\circ\sigma$ generated by
$1$ consists only of the numbers $1$ and $n-|\beta|$; as $\omega$
is connected, we must then have $n\leq 2$, which is not the case.
Therefore, the only possibilities for the pair $(|\omega\cap
m|,|\omega\cap b|)$ are the ones listed in the lemma.
\end{proof}

\begin{figure}
\Figw{3.6in}{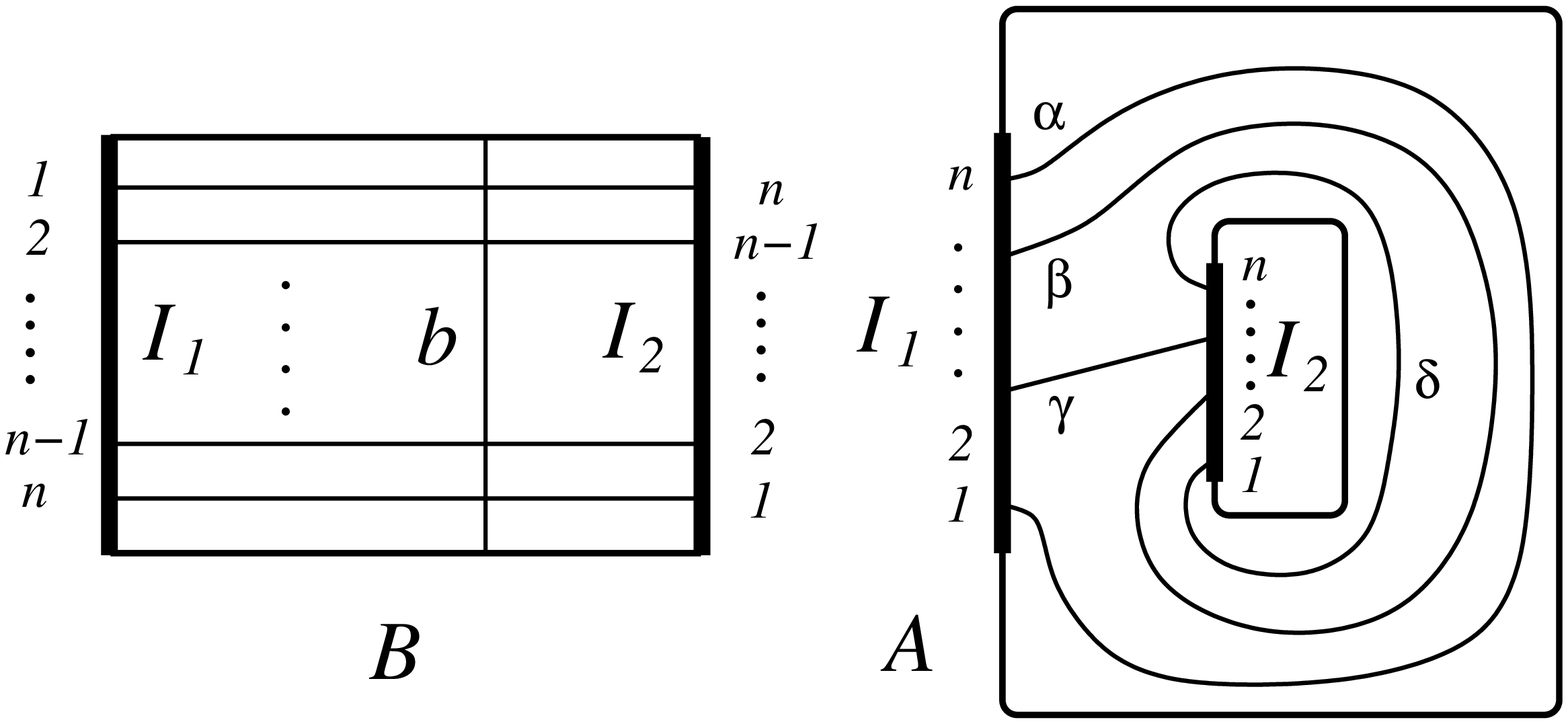}{The arcs $\omega\cap B\subset B$ and
$\omega\cap A\subset A$}{fig04}
\end{figure}

For a knot $L\subset S^3$, we will use the notation $C_2(L)$ to
generically denote any 2--cable of $L$; observe that any nontrivial
cable knot is prime.

\begin{lem}\label{centers}
For $1\leq s<n$, any center of $P_s$ is a knot of type
$$K_1\#\cdots\#K_s\#C_2(K_{s+1}\#\cdots\# K_n).$$
\end{lem}

\begin{proof}
By Lemma~\ref{circles}, any center $c_s\subset P_s$ can be
constructed as the union of two arcs: one that runs along the band
$B_n$ and the other any spanning arc of $A_s$. Since the annulus
`swallows' the factor $K_1\#\cdots\#K_s$ and `follows'
$K_{s+1}\#\cdots\# K_n$, any  such center circle $c_s$ (shown in
Figure~\ref{fig02}) isotopes into a knot of the type represented in
Figure~\ref{fig03}, which has the given form.
\end{proof}

\begin{figure}
\Figw{2.5in}{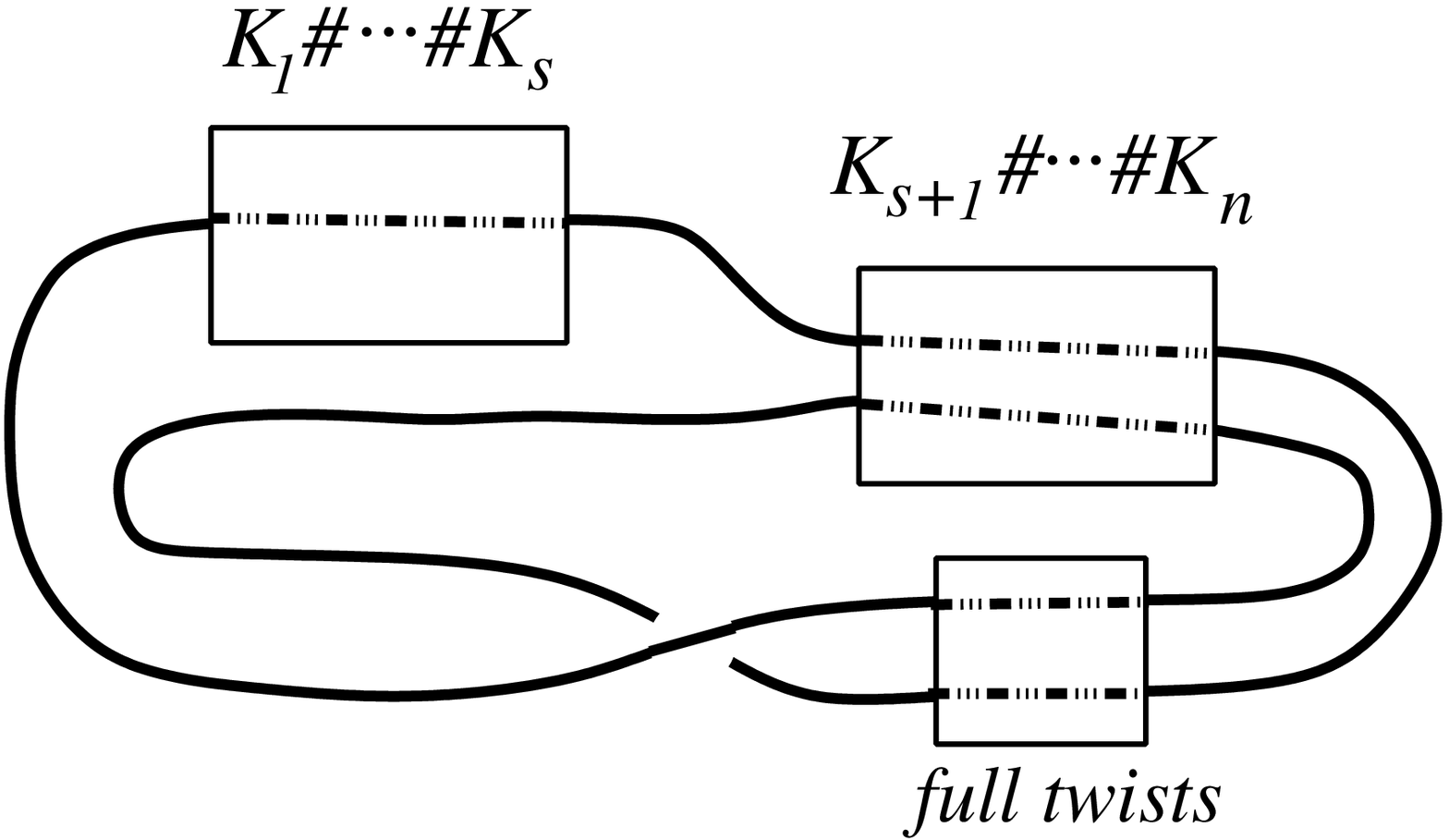}{The knot $c_s\subset S^3$}{fig03}
\end{figure}

\begin{lem}\label{notiso1}
For $1\leq r<s<n$, there is no homeomorphism $f\co S^3\rightarrow
S^3$ which maps $P_r$ onto $P_s$.
\end{lem}

\begin{proof}
Suppose there is a homeomorphism $f\co S^3\rightarrow S^3$ with
$f(P_r)=P_s$; then for any center $c_r$ of $P_r$, $c_s=f(c_r)$ is
a center of $P_s$. But then $c_r$ and $c_s$ have the same knot
type in $S^3$, which by Lemma~\ref{centers} can not be the case
since $c_r$ has $r+1$ prime factors while $c_s$ has $s+1$ prime
factors.
\end{proof}

\begin{proof}[Proof of Theorem~\ref{t3}]
By Lemma~\ref{notiso1}, the Seifert Klein bottles $P_r$ and $P_s$
for $\mc{K}_n$ are not equivalent for $1\leq r<s<n$, hence
$|\sk(\mc{K}_n,0)|\geq n-1$ and the theorem follows.
\end{proof}

\section{Cable and composite knots}\label{cable-comp}

In this section we assume that $K$ is a nontrivial knot in $S^3$
whose exterior $X_K$ contains an essential annulus $A$ and an
essential Seifert Klein bottle $P$; that is, $K$ is a crosscap
number two cable or composite knot. We assume that $A$ and $P$
have been isotoped so as to intersect transversely with $|A\cap
P|$ minimal, and denote by $G_A=A\cap P\subset A$ and $G_P=A\cap
P\subset P$ their graphs of intersection. We classify these graphs
in the next lemma; the case when $A$ has meridional boundary slope
is treated in full generality in \cite[Lemma 7.1]{eisner1}.

\begin{lem}\label{zerotwo}
Either $A\cap P=\emptyset$ or $\Delta(\partial A,\partial P)=1$
and $A\cap P$ consists of a single arc which is spanning in $A$
and separates $P$ into two Moebius bands.
\end{lem}

\begin{proof}
Suppose first that $\Delta(\partial A,\partial P)=0$, so that
$\partial A\cap\partial P=\emptyset$; in particular, the boundary
slopes of $A$ and $P$ are integral and $K$ is a cable knot with
cabling annulus $A$. If $A\cap P\neq\emptyset$ then $A\cap P$
consists of nontrivial orientation preserving circles in $A$ and
$P$. If any such circle $\gamma\subset A\cap P$ is parallel to
$\partial P$ in $P$, we may assume it cobounds an annulus
$A_{\gamma}$ with $\partial P$ in $P$ such that $A\cap\intr
A_{\gamma}=\emptyset$, so, by minimality of $|A\cap P|$,
$A_{\gamma}$ must be an essential annulus in the closure $V$ of
the component of $X_K\setminus A$ containing it; as $K$ is a cable
knot, $V$ must be the exterior of some nontrivial knot of whom $K$
is a $q$--cable for some $q\geq 2$. But then the boundary slope of
$A_{\gamma}$ in $V$ is of the form $a/q$, that is, nonintegral nor
$\infty$, contradicting the fact that $A_{\gamma}$ is essential.
Hence no component of $A\cap P$ is parallel to $\partial P$ in $P$
and so $\wh{P}$ compresses in $K(\partial P)$ along a subdisk of
$\wh{A}$, which is not possible by Lemma~\ref{pinc}. Therefore $A$
and $P$ are disjoint in this case.

Suppose now that $\Delta(\partial A,\partial P)\neq 0$; by
minimality of $|A\cap P|$, $A\cap P$ consists only of arcs which
are essential in both $A$ and $P$. If $\alpha$ is one such arc
then, as $\alpha$ is a spanning arc of $A$, $|\partial P|=1$, and
$X_K$ is orientable, it is not hard to see that $\alpha$ must be a
{\it positive} arc in $P$, in the sense of \cite{tera1} (this fact
does not follow directly from the parity rule in \cite{tera1}
since $|\partial A|\neq 1$, but its proof is equally direct). Thus
all the arcs of $G_P$ are positive in $P$.

Suppose $a,b$ are arcs of $G_P$ which are parallel and adjacent in
$P$, and let $R$ be the closure of the disk component of
$P\setminus(a\cup b)$. Then $R$ lies in the closure of some
component of $X_K\setminus A$ and, by minimality of $|A\cap P|$,
the algebraic intersection number $\partial R\cdot\text{core}(A)$
must be $\pm 2$. But then $K$ must be a 2--cable knot, with cabling
annulus $A$, contradicting Lemma~\ref{kleiness} since $P$ is
essential. Therefore no two arcs of $G_P$ are parallel.

Since any two disjoint positive arcs in $P$ are mutually parallel,
and $A$ separates $X_K$, the above arguments show that $G_P$
consists of exactly one essential arc which is separating in $P$.
The lemma follows.
\end{proof}

\begin{proof}[Proof of Theorem~\ref{t4}]
By Lemma~\ref{zerotwo}, $A\cap P=\emptyset$ or $|A\cap P|=1$.
Throughout the proof, none of the knots considered will  be  a
2--cable, hence any Seifert Klein bottle constructed for them will
be essential by Lemma~\ref{kleiness}. Let $V,W$ be the closures of
the components of $X_K\setminus A$.

\setcounter{case}{0}
\begin{case}
$A\cap P=\emptyset$
\end{case}

Here $K$ must be a cable knot: for the slope of $\partial A$ must
be integral or $\infty$, the latter case being impossible since
otherwise $\wh{P}$ is a closed Klein bottle in $K(\infty)=S^3$.
Hence $\partial V$ and $\partial W$ are parallel tori in $S^3$,
and we may regard them as identical for framing purposes. We will
assume $P\subset V$.

Suppose that $V$ is a solid torus; then $P$ is not
$\pi_1$--injective in $X_K$. Let $\mu,\lambda$ be a standard
meridian--longitude pair for $\partial V$, framed as the boundary
of the exterior of a core of $V$ in $S^3$. Since $V$ is a solid
torus and $P$ is essential in $X_K$, $P$ must boundary compress in
$V$ to some Moebius band $B$ in $V$ with $\Delta(\partial
P,\partial B)=2$. Suppose the slope of $\partial P$ in $V$ is
$p\mu+q\lambda$; as the slope of $\partial B$ is of the form
$(2m+1)\mu+2\lambda$ for some integer $m$, we must have
$(2m+1)q-2p=\pm 2$, hence $q\equiv 0\mod 4$ and $p$ is odd.
Therefore, the slope of $\partial P$ (and hence that of $\partial
A$) must be of the form $(2(2m+1)n\pm 1)\mu+4n\lambda$ for some
$n\neq 0$. Conversely, it is not hard to see that any such slope
bounds an incompressible once-punctured Klein bottle in $V$, which
is easily seen to be unique up to equivalence. Therefore, $K$ is a
$(2(2m+1)n\pm 1,4n)$ cable of the core of $V$ for some integers
$m,n,
\ n\neq 0$, and (a) holds.

If $V$ is not a solid torus then $W$ is a solid torus and, since
$\partial V$ and $\partial W$ are parallel in $S^3$, we can frame
$\partial V$ as the boundary of the exterior of a core of $W$ in
$S^3$ via a standard meridian--longitude pair $\mu,\lambda$. Then a
core of $A$ has slope $p\mu+q\lambda$ in $\partial V$ with
$|q|\geq 2$, so $P$ has nonintegral boundary slope in $V$ and
hence must boundary compress in $V$ by Lemmas~\ref{pinc} and
\ref{kleiness} to an essential Moebius band $B$ in $V$ with
$\Delta(\partial P,\partial B)=2$; in particular, $P$ is not
$\pi_1$--injective in $V$, hence neither in $X_K$. The slope of
$\partial B$ in $\partial V$ must be of the form
$2(2m+1)\mu+\lambda$ for some integer $m$, and so the slope of
$\partial P$ in $\partial V$ is be of the form $(2(2m+1)(2n+1)\pm
2)\mu+(2n+1)\lambda$ for some integer $n$. It follows that $K$ is
a $(2(2m+1)(2n+1)\pm 2,2n+1)$ cable of some $(2m+1,2)$ cable knot
for some integers $m,n$, and (b) holds.

\begin{case}
$|A\cap P|=1$
\end{case}

By Lemma~\ref{zerotwo}, $B_1=P\cap V$ and $B_2=P\cap W$ are
Moebius bands whose boundaries intersect $A$ in a single spanning
arc.

If the slope of $\partial A$ in $X_K$ is $\infty$ then $K=K_1\#
K_2$ for some nontrivial knots $K_1,K_2$ with $V=X_{K_1}$ and
$W=X_{K_2}$. Moreover, as $B_i\subset X_{K_i}$, the $K_i$'s are
2--cable knots and (d) holds.

Suppose now that $\partial A$ has integral slope in $X_K$; then at
least one of $V,W$, say $V$, is a solid torus. As in Case 1, we
may regard $\partial V$ and $\partial W$ as identical tori for
framing purposes. If, say,  $W$ is not a solid torus, frame
$\partial V$ and $\partial W$ via a standard meridian--longitude
pair $\mu,\lambda$ as the boundary of the exterior of a core of
$V$ in $S^3$. Then $\partial A$ has slope $p\mu+q\lambda$ for some
integers $p,q$ with $|q|\geq 2$ in $V$ and $W$. Also, $\partial
B_1$ has slope $r\mu+2\lambda$ in $V$ while $\partial B_2$ has
slope $2s\mu+\lambda$ in $W$, for some odd integers $r,s$. As
$\Delta(\partial A,\partial B_i)=1$ for $i=1,2$, we must have
$2p-rq =\varepsilon$ and $p-2sq =\delta$ for some
$\varepsilon,\delta=\pm 1$. Hence $p=2sq+\delta$ and so
$2(2sq+\delta)-rq=\varepsilon$, that is,
$|(4s-r)q|=|\varepsilon-2\delta|=1$ or $3$. As $|q|\geq 2$, it
follows that $|q|=3$, hence $p=6s\pm 1$. Therefore, $K$ is a
$(6s\pm 1,3)$ cable of a $(s,2)$ cable knot for some odd integer
$s$, and (c) holds.

The situation is somewhat different if both $V$ and $W$ are solid
tori. Here one can find meridian--longitude framings $\mu,\lambda$
for $\partial V$ and $\partial W$, such that $\mu$ bounds a disk
in $V$ and $\lambda$ bounds a disk in $W$. Now, the slope of
$\partial A$ in $\partial V$ is of the form $p\mu+q\lambda$ for
some integers $p,q$ with $|p|,|q|\geq 2$. The slope of $\partial
B_1$ in $V$ is of the form $a\mu+2\lambda$, while that of
$\partial B_2$ in $W$ is of the form $2\mu+b\lambda$. As
$\Delta(\partial A,\partial B_i)=1$ for $i=1,2$, it follows that
$2p-aq=\varepsilon$ and $2q-bp=\delta$ for some
$\varepsilon,\delta=\pm 1$. Assuming, without loss of generality,
that $|p|>|q|\geq 2$, the only solutions to the above equations
can be easily shown to be $(|p|,|q|)=(5,3)$ or $(7,3)$. Hence $K$
must be the torus knot $T(\pm 5,3)$ or $T(\pm 7,3)$. Notice that
this case fits in Case (c) with $n=-1,0$.

Thus in each of the above cases the knot $K$ admits an essential
Seifert Klein bottle, which is unique in Case (a). We have also
seen that in Cases (a) and (b) such a surface is never
$\pi_1$--injective.

In case (c), $K$ is a cable knot with cabling annulus $A$ such
that $X_K=V\cup_A W$, where $V$ is a solid torus and $W$ is the
exterior of some (possibly trivial) knot in $S^3$, and $B_1=P\cap
V$, $B_2=P\cap W$ are Moebius bands. Using our notation for
$H(P),T_P,K'$ relative to the surface $P$, as
$N(P)=N(B_1)\cup_{A\cap N(P)} N(B_2)$, we can see that
$$H(P)=(V\setminus\intr N(B_1))\cup_D (W\setminus\intr N(B_2))$$
where $D=\text{cl}(A\setminus N(P\cap A))\subset A$ is a disk (a
rectangle). Observe that
\ben
\item[(i)] $V\setminus\intr N(B_1)=A_1\times I$, where $A_1=A_1\times 0$ is
the frontier annulus of $N(B_1)$ in $V$,

\item[(ii)] the rectangle $D\subset A$ has one side along one
boundary component of the annulus $A'_1=\partial (A_1\times
I)\setminus\intr A_1\subset \partial (A_1\times I)$  and the
opposite side along the other boundary component, and

\item[(iii)] $K'\cap (A_1\times I)$ is an arc with one endpoint on
each of the sides of the rectangle $\partial D$ interior to
$A'_1$.
\een

Let $\alpha$ be a spanning arc of $A'_1$ which is parallel and
close to one of the arcs of $\partial D$ interior to $A'_1$, and
which is disjoint from $D$; such an arc exists by (i)--(iii), and
$\alpha$ intersects $K'$ transversely in one point by (iii).
Therefore $\alpha\times I\subset A_1\times I$ is a properly
embedded disk in $H(P)$ which intersects $K'$ transversely in one
point, so $T_P$ is boundary compressible in $X_K$ and hence $P$ is
not $\pi_1$--injective.

In case (d), keeping the same notation as above, if $P$ is not
$\pi_1$--injective then $T_P$ is not $\pi_1$--injective either, so
$T_P$ compresses in $X_K$, and in fact in $H(P)$, producing a
surface with at least one component an annulus $A_P\subset H(P)$
properly embedded in $X_K$ with the same boundary slope as $P$.
Notice that $A_P$ can not be essential in $X_K$: for $A_P$ must
separate $X_K$, and if it is essential then not both graphs of
intersection $A\cap A_P\subset A$, $A\cap A_P\subset A_P$ can
consist of only essential arcs by the Gordon--Luecke parity rule
\cite{cgls}. Therefore $A_P$ must be boundary parallel in $X_K$
and the region of parallelism must lie in $H(P)$ by
Lemma~\ref{lifts}(d), so $T_P$ is boundary compressible and there
is a disk $D'$ in $H(P)$ intersecting $K'$ transversely in one
point. As before,
$$H(P)=(V\setminus\intr N(B_1))\cup_D (W\setminus\intr N(B_2)),$$
where  $V=X_{K_1}$ and $W=X_{K_2}$. Isotope $D'$ so as to
intersect $D$ transversely with $|D\cap D'|$ minimal. If $|D\cap
D'|>0$ and $E'$ is an outermost disk component of $D'\setminus D$,
then, as $|D'\cap K'|=1$, $E'$ is a nontrivial disk in, say,
$V\setminus\intr N(B_1)$ with $|E'\cap K'|=0$ or $1$. Hence
$V\setminus\intr N(B_1)$ is a solid torus whose boundary
intersects $K'$ in a single arc $K'_1$ with endpoints on
$D\subset\partial V$, and using $E'$ it is possible to construct a
meridian disk $E''$ of $V\setminus\intr N(B_1)$ disjoint from $D$
which intersects  $K'$ coherently and transversely in one or two
points. If $|D\cap D'|=0$ we set $E''=D'$.

Let $L_1$ be the trivial knot whose exterior $X_{L_1}$ is the
solid torus $V\setminus\intr N(B_1)$, so that $K_1$ is a 2--cable
of $L_1$. As $X_{K_1}=X_{L_1}\cup N(B_1)$, where the gluing
annulus $X_{L_1}\cap N(B_1)$ is disjoint from the arc $K'\cap
X_{L_1}$ in $\partial X_{L_1}$, $E''$ must also intersect
$\partial B_1$ coherently and transversely in one or two points.
In the first case $K_1$ must be a trivial knot, while in the
second case $X_{K_1}\subset S^3$ contains a closed Klein bottle.
As neither option is possible, $P$ must be $\pi_1$--injective.
\end{proof}

\begin{proof}[Proof of Corollary~\ref{torus}]
That any crosscap number two torus knot is of the given form
follows from the proof of Theorem~\ref{t4}. The uniqueness of the
slope bounded by a Seifert Klein bottle in each case follows from
\cite{valdez6}, and that no such surface is $\pi_1$--injective
also follows from the proof of Theorem~\ref{t4}. For a knot $K$ of
the form $T(\cdot,4n)$, any Seifert Klein bottle $P$ is disjoint
from the cabling annulus and can be constructed on only one side
of the cabling annulus. Since $P$ is not $\pi_1$--injective, $T_P$
compresses in $H(P)$ giving rise to the cabling annulus of $K$;
thus uniqueness and unknottedness follows. For the knots $T(\pm
5,3),T(\pm 7,3)$ any Seifert Klein bottle $P$ is separated by the
cabling annulus into two Moebius bands; as a Moebius band in a
solid torus is unique up to ambient isotopy fixing its boundary,
uniqueness of $P$ again follows. Since, in the notation of the
proof of Theorem~\ref{t4}, $H(P)=(V\setminus\intr N(B_1))\cup_D
(W\setminus\intr N(B_2)),$ and $V\setminus\intr N(B_1), \
W\setminus\intr N(B_2)$ are solid tori, $H(P)$ is a handlebody and
so any Seifert Klein bottle in these last cases is unknotted.
\end{proof}

\section{Primitives, powers, and companion annuli}\label{handlebodies}

Let $M^3$ be a compact orientable 3--manifold with boundary, and
let $A$ be an annulus embedded in $\partial M^3$. We say that a
separating annulus $A'$ properly embedded in $M^3$ is a {\it
companion} of $A$ if $\partial A'=\partial A$ and $A'$ is not
parallel into $\partial M^3$; we also say that $A'$ is a companion
of any circle $c$ embedded in $\partial M^3$ which is isotopic to
a core of $A$. Notice that the requirement of a companion annulus
being separating is automatically met whenever $M^3\subset S^3$,
and that if $\partial M^3$ has no torus component then we only
have to check that $A'$ is not parallel into $A$. The following
general result will be useful in the sequel.

\begin{lem}\label{soltor}
Let  $M^3$ be an irreducible and atoroidal 3--manifold with
connected boundary, and let $A',B'$ be companion annuli of some
annuli $A,B\subset\partial M^3$. Let $R,S$ be the regions in $M^3$
cobounded by $A,A'$ and $B,B'$, respectively. If  $A$ is
incompressible in $M^3$ then $R$ is a solid torus, and if $A$ and
$B$ are isotopic in $\partial M^3$ then $R\cap S\neq\emptyset$. In
particular, $A'$ is unique in $M^3$ up to isotopy.
\end{lem}

\begin{proof}
Let $c$ be a core of $A$; push $R$ slightly into $\intr M^3$ via a
small collar $\partial M^3\times I$ of $\partial M^3=\partial
M^3\times 0$, and let $A''=c\times I$. Observe the annulus $A''$
has its boundary components $c\times 0$ on $\partial M^3$ and
$c\times 1$ on $\partial R$. Also, $\partial R$ can not be
parallel into $\partial M^3$, for otherwise $A'$ would be parallel
into $\partial M^3$; $M^3$ being atoroidal,  $\partial R$ must
compress in $M^3$.

Let $D$ be a nontrivial compression disk of $\partial R$ in $M^3$.
If $D$ lies in $N^3=M^3\setminus\intr R$ then $\partial R$
compresses in $N^3\cup(\text{2--handle along }c\times 0)$ along the
circles $\partial D$ and $c\times 1\subset\partial R$, hence
$\partial D$ and $c\times 1$ are isotopic in $\partial R$, which
implies that $A$ compresses in $M^3$, contradicting our
hypothesis. Thus $D$ lies in $R$, so $R$ is a solid torus.

Suppose now that $R\cap S=\emptyset$, so that $A\cap
B=\emptyset=A'\cap B'$, and that $A,B$ are isotopic in $\partial
M^3$. Then one boundary component of $A$ and one of $B$ cobound an
annulus $A^*$ in $\partial M^3$ with interior disjoint from $A\cup
B$ (see Figure~\ref{lem51}(a)). Since, by the above argument, $R$
and $S$ are solid tori with $A,B$ running more than once around
$R,S$, respectively, the region $N(R\cup S\cup A^*)$ is a Seifert
fibered space over a disk with two singular fibers, hence not a
solid torus, contradicting our initial argument. Therefore $R\cap
S\neq\emptyset$.

Finally, isotope $B$ into the interior of $A$, carrying $B'$ along
the way so that $A'$ and $B'$ intersect transversely and
minimally. If $A'\cap B'=\emptyset$ then $B'$ lies in the solid
torus $R$ and so must be parallel to $A'$ in $R$. Otherwise, any
component $B''$ of $B'\cap (M^3\setminus\intr R)$ is an annulus
that cobounds a region $V\subset M^3\setminus\intr R$ with some
subannulus $A''$ of $A'$ (see Figure~\ref{lem51}(b)). Since $|A'\cap
B'|$ is minimal, $A''$ and $B''$ are not parallel within $V$ and
so $R\cup_{A''} V$ is not a solid torus. But then the frontier of
$R\cup_{A''} V$ is a companion annulus of $A$, contradicting our
initial argument. The lemma follows.
\end{proof}

\begin{figure}\nocolon
\Figw{3.3in}{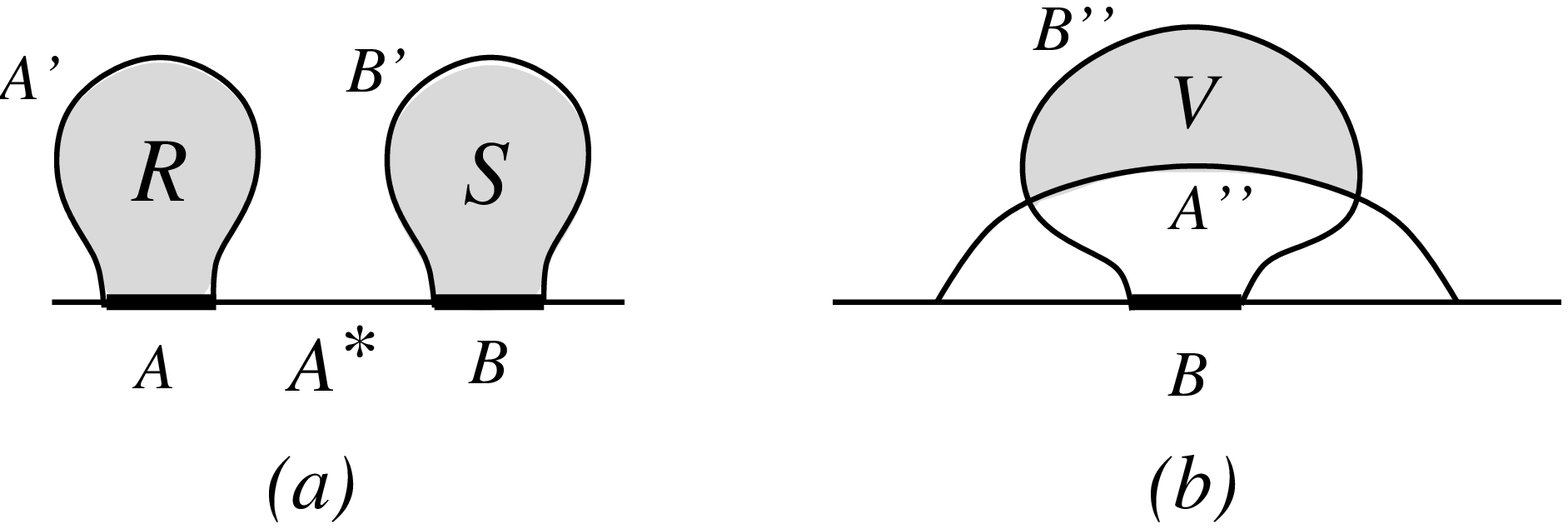}{}{lem51}
\end{figure}

\begin{rem}
In the context of Lemma~\ref{soltor}, if $A\subset\partial M^3$
compresses in $M^3$ and $M^3$ is irreducible, then $A$ has a
companion iff the core of $A$ bounds a nonseparating disk in
$M^3$, in which case $A$ has infinitely many nonisotopic companion
annuli.
\end{rem}

In the special case when $H$ is a genus two handlebody, an
algebraic characterization of circles in $\partial H$ that admit
companion annuli will be useful in the sequel, particularly in
Section~\ref{pretzels}; we introduce some terminology in this
regard. For $H$ a handlebody and  $c$ a circle embedded in
$\partial H$, we say $c$ is {\it algebraically primitive} if $c$
represents a primitive element in $\pi_1(H)$ (relative to some
basepoint), and we say $c$ is {\it geometrically primitive} if
there is a disk $D$ properly embedded in $H$ which intersects $c$
transversely in one point. It is well known that these two notions
of primitivity for circles in $\partial H$ coincide, so we will
refer to such a circle $c$ as being simply {\it primitive} in $H$.
We say that $c$ is a {\it power} in $H$ if $c$ represents a proper
power of some nontrivial element of $\pi_1(H)$.

The next result follows essentially from \cite[Theorem
4.1]{cassongor}; we include a short version of the argument for
the convenience of the reader.

\begin{lem}\label{comp}
Let $H$ be a genus two handlebody and $c$ a circle embedded in
$\partial H$ which is nontrivial in $H$. Then,
\ben
\item[\rm(a)] $\partial H\setminus c$ compresses in $H$ iff $c$ is
primitive or a proper power in $H$, and

\item[\rm(b)] $c$ has a companion annulus in $H$ iff $c$ is a
power in $H$.
\een
\end{lem}

\begin{proof}
For (a), let $D\subset H$ be a compression disk of $\partial
H\setminus c$. If $D$ does not separate $H$ then, since $c$ and
$\partial D$ can not be parallel in $\partial H$, there is a
circle $\alpha$ embedded in $\partial H\setminus c$ which
intersects $\partial D$ transversely in one point, and so the
frontier of a regular neighborhood in $H$ of $D\cup\alpha$ is a
separating compression  disk of $\partial H\setminus c$. Thus, we
may assume that $D$ separates $H$ into two solid tori $V,W$ with
$D=V\cap W=\partial V\cap\partial W$ and, say, $c\subset V$.
Therefore $c$ is either primitive or a power in $V$ and hence in
$H$. Conversely, suppose $c$ is primitive or a power in $H$. Then
$\pi_1(H(c))$ is either free cyclic or has nontrivial torsion by
\cite[Theorems N3 and 4.12]{mks} and so the pair
$(H(c),\partial H(c))$ is not irreducible by
\cite[Theorem 9.8]{hempel}. Therefore, by the 2--handle addition
theorem (cf \cite{cassongor}), the surface $\partial H\setminus
c$ must be compressible in $H$.

For (b), let $A$ be an annulus neighborhood of $c$ in $\partial
H$. If $A'$ is a companion annulus of $A$ then $A'$ must boundary
compress in $H$ into a compression disk for $\partial H\setminus
A$; as in (a), we can assume that $\partial H\setminus A$
compresses along a disk $D$ which separates $H$ into two solid
tori $V,W$ with $D=V\cap W=\partial V\cap\partial W$ and $c,
A,A'\subset V$. Since $A'$ is parallel in $V$ into $\partial V$
but not into $A$, it follows that $c$ is a power in $V$, hence in
$H$. The converse follows in a similar way.
\end{proof}

A special family of incompressible pairs of pants properly
embedded in a 3--manifold $H\subset S^3$ with $\partial H$ a genus
two surface, which are not parallel into $\partial H$, appear
naturally in Section~\ref{hypknots}. We will establish some of
their properties in the next lemma; the following construction
will be useful in this regard. If $F$ is a proper subsurface of
$\partial H$, $c$ is a component of $\partial F$, and $A$ is a
companion annulus of $c$ in $H$ with $\partial A=\partial_1
A\cup\partial_2 A$, we isotope $A$ so that, say, $\partial_1 A=c$
and $\partial_2 A\subset\partial H\setminus F$, and denote by
$F\oplus A$ the surface $F\cup A$, isotoped slightly so as to lie
properly embedded in $H$.

For $c_1,c_2,c_3$ disjoint circles embedded in $\partial H$ and
nontrivial in $H$, we say $c_1,c_2$ are {\it simultaneously
primitive away from $c_3$} if there is some disk $D$ in $H$
disjoint from $c_3$ which transversely intersects $c_1$ and $c_2$
each in one point; notice that if $H$ is a handlebody and
$c_1,c_2$ are simultaneously primitive away from some other
circle, then $c_1,c_2$ are indeed primitive in $H$. We also say
$c_1,c_2$ are {\it coannular} if they cobound an annulus in $H$.

\begin{lem}\label{pants}
Let $H\subset S^3$ be a connected atoroidal 3--manifold with
connected boundary of genus two. Let $c_1,c_2,c_3$ be disjoint
nonseparating circles embedded in $\partial H$ which are
nontrivial in $H$, no two are coannular, and separate $\partial H$
into two pairs of pants $P_1,P_2$. Let $Q_0$ be a pair of pants
properly embedded in $H$ with $\partial Q_0=c_1\cup c_2\cup c_3$
and not parallel into $\partial H$. Then $Q_0$ is incompressible
and separates $H$ into two components with closures $H_1,H_2$ and
$\partial H_i=Q_0\cup P_i$, and if $Q_0$ boundary compresses in
$H$ the following hold:
\ben
\item[\rm(a)] $c_i$ has a companion annulus in $H$ for some $i=1,2,3$;

\item[\rm(b)] if only $c_1$ has a companion annulus in $H$, say
$A'_1$, then,
\ben
\item[\rm(i)] for $i,j=1,2$, $P_i\oplus A'_1$ is isotopic to
$P_j$ in $H$ iff $H$ is a handlebody and $c_2,c_3$ are
simultaneously primitive in $H$ away from $c_1$;

\item[\rm(ii)]
if $Q_0$ boundary compresses in $H_i$ then it boundary compresses
into a companion annulus of $c_1$ in $H_i$, $Q_0=P_i\oplus A'_1$
in $H$, $H_i$ is a handlebody, and any pair of pants in $H_i$ with
boundary $\partial Q_0$ is parallel into $\partial H_i$,

\item[\rm(iii)] if  $R_0$ is a boundary compressible pair of
pants in $H$ disjoint from $Q_0$ with $\partial R_0$ isotopic to
$\partial Q_0$ in $\partial H$, then $R_0$ is parallel to $Q_0$ or
$\partial H$.
\een
\een
\end{lem}

\begin{proof}
Observe $H$ is irreducible and $Q_0$ is incompressible in $H$;
since $H\subset S^3$, $Q_0$ must separate $H$, otherwise $Q_0\cup
P_i$ is a closed nonseparating surface in $S^3$, which is
impossible. Also, by Lemma~\ref{soltor}, a companion annulus of
any $c_i$ is unique up to isotopy and cobounds a solid torus with
$\partial H$. Let $H_1,H_2$ be the closures of the components of
$H\setminus Q_0$, with $\partial H_i=Q_0\cup P_i$; as $Q_0$ is
incompressible, both $H_1$ and $H_2$ are irreducible and
atoroidal, so again companion annuli in $H_i$ are unique up to
isotopy and cobound solid tori with $\partial H_i$. Let $D$ be a
boundary compression disk for $Q_0$, say $D\subset H_1$. We
consider three cases.

\setcounter{case}{0}
\begin{case}\label{cas1}
The arc $Q_0\cap\partial D$ does not separate $Q_0$.
\end{case}

Then $D$ is a nonseparating disk in $H_1$, and we may assume the
arc $Q_0\cap\partial D$ has one endpoint in $c_1$ and the other in
$c_2$, so that $|c_1\cap D|=1=|c_2\cap D|$ and $c_3\cap
D=\emptyset$. Hence the frontier $D'$ of a small regular
neighborhood of $c_1\cup D$ is a properly embedded separating disk
in $H_1$ which intersects $c_2$ in two points and whose boundary
separates $c_1$ from $c_3$ in $\partial H$; the situation is
represented in Figure~\ref{lem53ab}, with $c_1=u,c_2=v,c_3=w$ and
$D=\Sigma,D'=\Sigma'$. Clearly, boundary compressing $Q_0$ along
$D$ produces an annulus $A'_3$ in $H_1$ with boundaries parallel
to $c_3$ in $\partial H_1$, and if $A'_3$ is parallel into
$\partial H_1$ then $Q_0$ itself must be parallel into
$P_1\subset\partial H_1$, which is not the case. Therefore $c_3$
has a companion annulus in $H_1$, hence in $H$.

Notice that $H_1$ must be a handlebody in this case since the
region cobounded by $A_3'$ and $\partial H_1$ is a solid torus by
Lemma~\ref{soltor}, and that $Q_0=P_1\oplus A_3'$ in $H$ and
$c_1,c_2$ are simultaneously primitive in $H_1$ away from $c_3$.

\begin{figure}\nocolon
\Figw{2.2in}{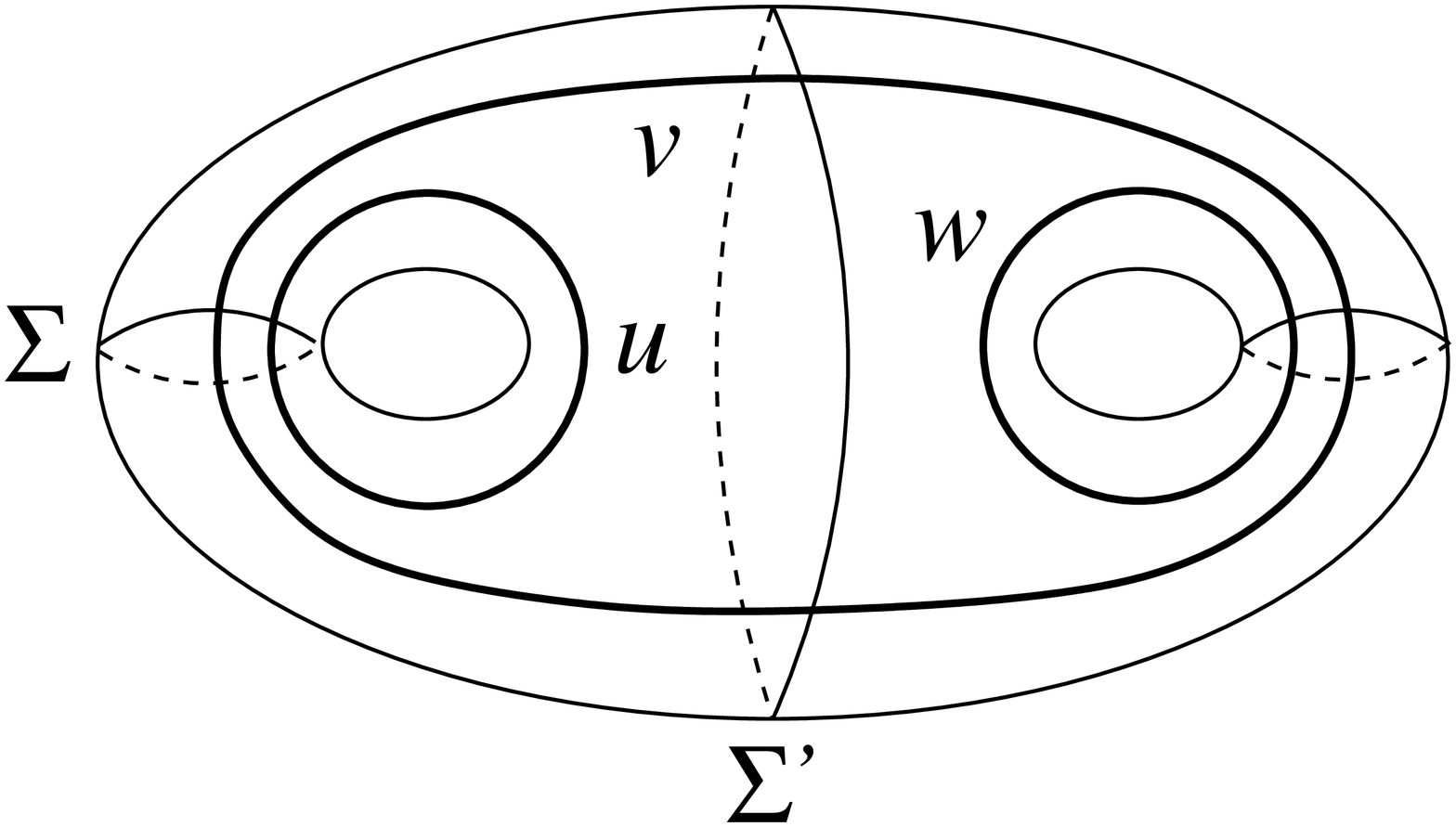}{}{lem53ab}
\end{figure}

\begin{case}\label{cas2}
The arc $Q_0\cap\partial D$ separates $Q_0$ and $D$ separates
$H_1$.
\end{case}

Here the endpoints of the arc $Q_0\cap\partial D$ lie in the same
component of $\partial Q_0$, so we may assume that $|c_1\cap D|=2$
with  $c_1\cdot D=0$ while $c_2,c_3$ are disjoint from and
separated by $D$; the situation is represented in
Figure~\ref{lem53ab}, with $c_1=v,c_2=u,c_3=w$, and $D=\Sigma'$.
Thus, boundary compressing $Q_0$ along $D$ produces two annuli
$A'_2,A'_3$ in $H_1\setminus D$ with boundaries parallel to
$c_2,c_3$, respectively. Since $Q_0$ is not parallel into
$P_1\subset\partial H_1$, at least one of these annuli must be a
companion annulus.

Notice that if only one such annulus, say $A'_2$, is a companion
annulus then, as in Case 1, $H_1$ is a handlebody, $Q_0=P_1\oplus
A'_2$ in $H$,  and $c_1,c_3$ are simultaneously primitive in $H_1$
away from $c_2$.

\begin{case}
The arc $Q_0\cap\partial D$ separates $Q_0$ and $D$ does not
separate $H_1$.
\end{case}

As in Case~\ref{cas2}, we may assume that $|c_1\cap D|=2$  with
$c_1\cdot D=0$ while $c_2,c_3$ are disjoint from $D$. Since $D$
does not separate $H_1$, boundary compressing $Q_0$ along $D$
produces two nonseparating annuli in $H_1$, each with one boundary
parallel to $c_2$ and the other parallel to $c_3$. Since $c_2,c_3$
are not coannular in $H$,  this case does not arise. Therefore (a)
holds.

For (b)(i), let $V,H'$ be the closures of the  components of
$H\setminus A'_1$, with $V$ a solid torus and $P_j\subset \partial
H'$; observe that $\partial H'$ can be viewed as $(P_i\oplus
A'_1)\cup P_j$. If $P_i\oplus A'_1$ and $P_j$ are isotopic in $H$
then $H'\approx P_j\times I$ with $P_j$ corresponding to
$P_j\times 0$, from which it follows that $c_2,c_3$ are
simultaneously primitive in $H$ away from $c_1$. Moreover, $c_1$
is primitive in the handlebody $H'$, and so $H=H'\cup_{A'_1}V$ is
also a handlebody. Conversely, suppose $H$ is a handlebody and
$c_2,c_3$ are simultaneously primitive away from $c_1$; notice
that $c_2,c_3,A'_1\subset\partial H'$. Suppose $D$ is a disk in
$H$ realizing the simultaneous primitivity of $c_2,c_3$ away from
$c_1$; we assume, as we may, that $D$ lies in $H'$ (see
Figure~\ref{lem53ab} with $D=\Sigma, \ c_2=u, \ c_3=v, \
\text{core}(A'_1)=w$). Compressing $\partial H'\setminus\intr A'_1$
in $H'$ along $D$ gives rise to an annulus $A''_1$ in $H$ which,
due to the presence of $A'_1$, is necessarily a companion annulus
of $c_1$ in $H$. It follows from Lemma~\ref{soltor} that $A''_1$
and $A'_1$ are parallel in $H'$, hence that $P_i\oplus A'_1$ and
$P_j$ are parallel in $H'$ and so in $H$.

For (b)(ii), if $Q_0$ boundary compresses into, say, $H_1$, then
it follows immediately from the proof of (a) that $H_1$ is a
handlebody, $Q_0=P_1\oplus A_1'$, and $c_2,c_3$ are simultaneously
primitive in $H_1$ away from $c_1$, so $Q_0$ compresses into a
companion annulus of $c_1$ in $H_1$. If $R_0$ is any pair of pants
in $H_1$ with boundary $\partial Q_0$ then by the same argument
$R_0$ must be isotopic in $H_1$ to $Q_0\oplus A'_1$ or $P_1\oplus
A'_1$, hence by (b)(i) $R_0$ is parallel into $\partial H_1$.

For (b)(iii), consider the disjoint pairs of pants $Q_0,R_0$ and
suppose $R_0$ is also not parallel into $\partial H$. By (b)(ii),
$Q_0=P_i\oplus A'_1$ and $R_0=P_j\oplus A'_1$ for some
$i,j\in\set{1,2}$. Thus, if $i\neq j$, $Q_0$ boundary compresses
in the direction of $P_i$, away from $R_0$, into a companion
annulus of $c_1$; a similar statement holds for $R_0$, and so such
companion annuli of $c_1$ are separated by $Q_0\cup R_0$,
contradicting Lemma~\ref{soltor}. Therefore $i=j$, so $Q_0$ is
parallel to $R_0$. The lemma follows.
\end{proof}

\begin{rem}
If any two of the circles $c_1,c_2,c_3$ in Lemma~\ref{pants} are
coannular in $H$ then part (a) need not hold; an example of this
situation can be constructed as follows. Let $H$ be a genus two
handlebody and let $c_1,c_2$ be disjoint, nonparallel,
nonseparating circles in $\partial H$ which are nontrivial in $H$
and cobound an annulus $A$ in $H$. Let $\alpha\subset\partial H$
be an arc with one endpoint in $c_1$ and the other in $c_2$ which
is otherwise disjoint from $A$. We then take the pair of pants
$Q_0$ to be the frontier of $H_1=N(A\cup\alpha)$ in $H$ so that,
up to isotopy, $\partial Q_0=c_1\cup c_2\cup K$, where
$K\subset\partial H$ is the sum of $c_1$ and $c_2$ along $\alpha$.
Observe $Q_0$ boundary compresses in $H_1$ as in Case 3 of
Lemma~\ref{pants}(a). Thus, if $c_1$ and $c_2$ are primitive in
$H$ then, by Lemma~\ref{comp}, no component of $\partial Q_0$ has
a companion in $H$ and $Q_0$ need not be parallel into $\partial
H$, as illustrated by the example in Figure~\ref{example52}.
Moreover, if $P$ is the once-punctured Klein bottle $A\cup B$,
where $B\subset\partial H$ is a band with core $\alpha$, pushed
slightly off $\partial H$ so as to properly embed in $H_1$, then
$\partial P$ is isotopic to $K$, $H_1$ is a regular neighborhood
of  $P$, and the two components of $\partial Q_0$  isotopic to
$c_1,c_2$ are lifts of the meridian of $P$.
\end{rem}

\begin{figure}
\Figw{3.8in}{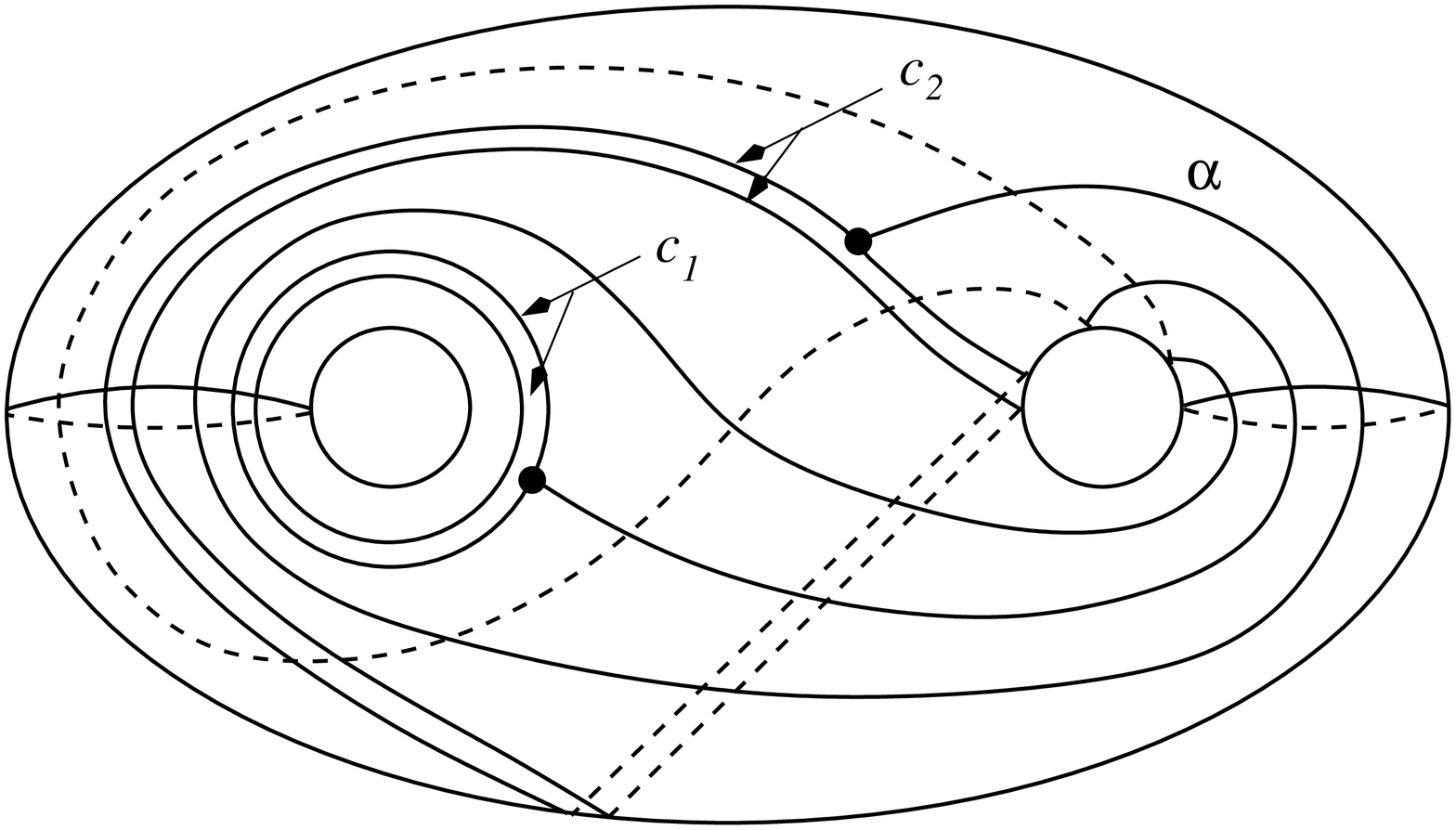}{Coannular primitive circles $c_1,c_2$ in the
handlebody $H$}{example52}
\end{figure}

\section{Hyperbolic knots}\label{hypknots}

In this section we fix our notation and let $K$ be a hyperbolic
knot in $S^3$ with exterior $X_K$. If $P,Q$ are distinct elements
of $\sk(K,r)$ which have been isotoped so as to intersect
transversely and minimally, then $|P\cap Q|>0$, $\partial
P\cap\partial Q=\emptyset$, and each circle component of $P\cap Q$
is nontrivial in $P$ and $Q$. Notice that any circle component
$\gamma$ of $P\cap Q$ must be orientation preserving in both $P,
Q$, or orientation reversing in both $P,Q$. If $\gamma$ is a
meridian (longitude, center) in both $P$ and $Q$, we will say that
$\gamma$ is a {\it simultaneous} meridian (longitude, center,
respectively) in $P,Q$. Recall our notation for $H(P),T_P,K',A'_K$
relative to the surface $P$.

\begin{lem}\label{nocomp}
Let $m_1,m_2\subset T_P$ and $l\subset T_P$ be the lifts of a
meridian circle and a center circle of $P$, respectively. Then
$m_1,m_2$ can not both have companions in $H(P)$, and neither $l$
nor $K'\subset\partial H(P)$ have companions.
\end{lem}

\begin{proof}
Let $A_1,A_2,A$ be annular neighborhoods of $m_1,m_2,l$ in $T_P$,
respectively, with $A_1\cap A_2=\emptyset$. Let $A_1',A_2'$ be
companions of $A_1,A_2$ in $H(P)$, respectively, and suppose they
intersect transversely and minimally; then $A_1'\cap
A_2'=\emptyset$ by Lemma~\ref{lifts}(b). Let $V_i$ be the region
in $H(P)$ cobounded by $A_i,A'_i$ for $i=1,2$. Now, the lifts
$m_1,m_2$ cobound the annulus $A(m)$ in $N(P)$ for some meridian
circle $m\subset P$. Since $A'_i$ and $A_i$ are not parallel in
$V_i$, it follows that, for a small regular neighborhood
$N=A(m)\times I$ of $A(m)$ in $N(P)$, $V_1\cup N\cup V_2\subset
X_K$ is the exterior $X_L$ of some nontrivial knot $L$ in $S^3$
with $A_1\subset X_L$ an essential annulus (see
Figure~\ref{lifts2}(a)). Observe $H(P)$ is irreducible and atoroidal
since $K$ is hyperbolic and $P$ is essential. As $A_1,A_2$ are
incompressible in $H(P)$, the $V_i$'s are solid tori by
Lemma~\ref{soltor} and so $L$ is a nontrivial torus knot with
cabling annulus $A_1$. However, as the pair of pants
$P_0=P\setminus\intr X_L$ has two boundary components coherently
oriented in $\partial X_L$ with the same slope as $\partial A_1$,
$K$ must be a satellite of $L$ of winding number two,
contradicting the hyperbolicity of $K$.

\begin{figure}\nocolon
\Figw{3in}{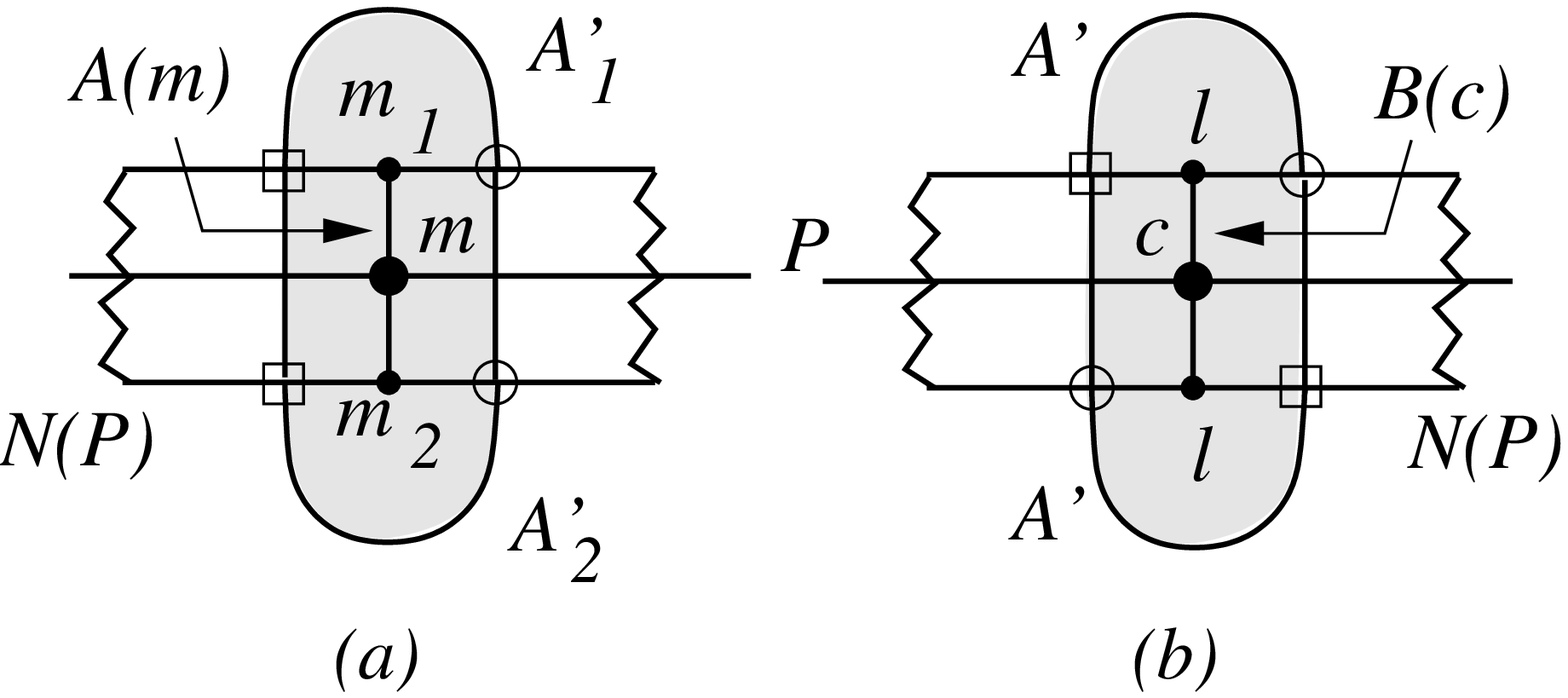}{}{lifts2}
\end{figure}

Suppose $A'$ is a companion of $A$ in $H(P)$, and let $V$ be the
region in $H(P)$ cobounded by $A,A'$. The circle $l$ bounds a
Moebius band $B(c)$ in $N(P)$ with $B(c)\cap P$ some center circle
$c$ of $P$. If $M$ is a small regular neighborhood of $B(c)$ in
$N(P)$ then, as $A,A'$ are not parallel in $V$, $V\cup M$ is the
exterior $X_L$ of some nontrivial knot $L$ in $S^3$ (see
Figure~\ref{lifts2}(b)), and $A$ is an essential annulus in $X_L$
which, since $M$ is a solid torus, necessarily has integral
boundary slope in $\partial X_L$. This time, $P_0=P\setminus\intr
X_L$ is a once-punctured Moebius band with one boundary component
in $\partial X_L$ having the same slope as $\partial A$, so $K$ is
a nontrivial satellite of $L$ with odd winding number, again
contradicting the hyperbolicity of $K$.

Finally, if $A'_K$ has a companion annulus $B'$ then, as $K$ is
hyperbolic, $B'$ must be boundary parallel in $X_K$ in the
direction of $P$, contradicting Lemma~\ref{lifts}(d). The lemma
follows.
\end{proof}

Given distinct elements $P,Q\subset\sk(K,r)$, if $P\cap Q$ is a
single simultaneous meridian or some pair of disjoint simultaneous
centers, we will say that $P$ and $Q$ intersect {\it meridionally}
or {\it centrally}, respectively.

\begin{lem}\label{mlc}
If $P,Q\in\sk(K,r)$  are distinct elements which intersect
transversely and minimally, then $P,Q$ intersect meridionally or
centrally.
\end{lem}

\begin{proof}
By minimality of $|P\cap Q|$, any component $\gamma$ of $P\cap Q$
must be nontrivial in both $P$ and $Q$, hence, in $P$ or $Q$,
$\gamma$ is either a circle parallel to the boundary, a meridian,
a longitude, or  a center circle.  Observe that if $\gamma$ is a
center in $P$ then, as it is orientation reversing in $P$ it must
be orientation reversing in $Q$, and hence $\gamma$ must also be a
center in $Q$.

Suppose $\gamma$ is parallel to $\partial P$ in $P$; without loss
of generality, we may assume that $\partial P$ and $\gamma$
cobound an annulus $A_P$ in $P$ with $Q\cap\intr A_P=\emptyset$.
Since $\gamma$ preserves orientation in $P$ it must also preserve
orientation in $Q$, hence $\gamma$ is either a meridian or
longitude of $Q$, or parallel to $\partial Q$. In the first two
cases, $\wh{Q}$ would compress in $K(\partial Q)$ via the disk
$\wh{A}_P$, which is not the case by Lemma~\ref{pinc}; thus,
$\gamma$ and $\partial Q$ cobound an annulus $A_Q$ in $Q$. As $K$
is hyperbolic, the annulus $A_P\cup_{\gamma} A_Q\subset H(P)$ is
boundary parallel in $X_K$ by Lemma~\ref{nocomp}; but then $|P\cap
Q|$ is not minimal, which is not the case. Therefore no component
of $P\cap Q$ is parallel to $\partial P,\partial Q$ in $P,Q$,
respectively.

Suppose now that $\gamma$ is a meridian in $P$ and a longitude in
$Q$. Then $P\cap Q$ consists only of meridians of $P$ and
longitudes of $Q$. If $\gamma$ is the only component of $P\cap Q$
then $P\asymp Q$ is a nonorientable (connected) surface properly
embedded in $X_K$ with two boundary components, which is
impossible. Thus we must have $|P\cap Q|\geq 2$; as the circles
$P\cap Q$ are mutually parallel meridians in $P$, it follows that
$P\cap H(Q)$ consists of a pair of pants and at least one annulus
component $A$. But $P\cap N(Q)$ consists of a disjoint collection
of annuli $\{A_i\}$ with $\{A_i\cap Q\}$ disjoint longitude
circles of $Q$, and so the circles $\cup\partial A_i$ form at most
two parallelism classes in $T_Q\subset\partial H(Q)$,
corresponding to the lifts of some disjoint pair of centers of
$Q$. Since the circles $\partial A$ are among those in
$\cup\partial A_i$, and $A$ is not parallel into $\partial H(Q)$
by minimality of $|P\cap Q|$, we contradict Lemmas~\ref{lifts} and
\ref{nocomp}.

Therefore, each component of $P\cap Q$ is a simultaneous meridian,
longitude, or center of $P,Q$. There are now two cases left to
consider.

\setcounter{case}{0}
\begin{case}
$P\cap Q$ consists of simultaneous meridians.
\end{case}
Suppose $|P\cap Q|=k+1$, $k\geq 0$. Then $Q\cap N(P)$ consists of
disjoint parallel annuli $A_0,\dots, A_k$, each intersecting $P$
in a meridian circle, and $Q\cap H(P)=Q_0\cup A'_1\cup\dots\cup
A'_k$, where $Q_0$ is a pair of pants with $\partial
Q\subset\partial Q_0$ and the $A'_i$'s are annuli, none of which
is parallel into $\partial H(P)$. The circles $\cup\partial
A_i=\cup\partial A_i'$ consist of two parallelism classes in
$\partial H$, denoted I and II, corresponding to the two distinct
lifts of a meridian circle of $P$ to $\partial N(P)$.

By Lemma~\ref{lifts}, the circles $\partial A_i'$ are both of type
I or both of type II, for each $i$. Also, the components of
$\partial Q_0$ are $\partial Q$ and two circles
$\partial_1Q_0,\partial_2Q_0$ of type I or II. If the circles
$\partial_1Q_0,\partial_2Q_0$ are both of type I or both of type
II, then the union of $Q_0$ and an annulus in $T_P$ cobounded by
$\partial_1Q_0$ and $\partial_2Q_0$ is a once-punctured surface in
$X_K$ disjoint from $P$, contradicting Lemma~\ref{lifts}(b).
Therefore, one of the circles $\partial_1Q_0,\partial_2Q_0$ is of
type I and the other of type II. Thus, if $k>0$ then some annulus
$A'_i$ has boundaries of type I and some annulus $A'_j$ has
boundaries of type II, which contradicts Lemma~\ref{nocomp}.
Therefore $k=0$, and so $P\cap Q$ consists of a single
simultaneous meridian.

\begin{case}
$P\cap Q$ consists of simultaneous longitudes or centers.
\end{case}

$P\cap Q$ can have at most two simultaneous centers; if it has at
most one simultaneous center then $P\asymp Q$ is a nonorientable
surface properly embedded in $X_K$ with two boundary components,
which is not possible. Therefore $P\cap Q$ contains a pair
$c_1,c_2$ of disjoint center circles, so $Q\cap N(P)$ consists of
two Moebius bands and, perhaps, some annuli, while $Q\cap H(P)$
consists of one pair of pants $Q_0$ and, perhaps, some annuli.
Thus the circles $Q\cap\partial H(P)$ are divided into two
parallelism classes, corresponding to the lifts of $c_1$ and
$c_2$, and  we may proceed as in Case 1 to show that $Q\cap H(P)$
has no annulus components. Hence $P\cap Q=c_1\cup c_2$.
\end{proof}

\begin{rem}
Notice that if $P,Q$ are any two distinct elements of $\sk(K,r)$,
so $P\cap Q$ is central or meridional, and  $Q_0=Q\cap H(P)$, then
by Lemma~\ref{lifts}, since $K$ is hyperbolic, $H(P)$ and
$\partial Q_0$ satisfy the hypothesis of Lemma~\ref{pants};
however, $Q_0$ need not boundary compress in $H(P)$.
\end{rem}

The following result gives constraints on the exteriors of distinct
elements of $\sk(K,r)$.

\begin{lem}\label{hbdy-irred}
Suppose $P,Q\in\sk(K,r)$ are distinct elements which intersect
centrally or meridionally; let $Q_0=Q\cap H(P)$, and let $V,W$ be
the closures of the components of $H(P)\setminus Q_0$. Then $Q_0$
is not parallel into $\partial H(P)$, and
\ben
\item[\rm(a)] for $X=V,W,$ or $H(P)$, either $X$ is a handlebody or
the pair $(X,\partial X)$ is irreducible and atoroidal;

\item[\rm(b)] if $P\cap Q$ is central then
$(H(P),\partial H(P))$ is irreducible and $Q_0$ is boundary
incompressible in $H(P)$, and

\item[\rm(c)] if the pair $(V,\partial V)$ is irreducible then
$W$ is a handlebody.
\een
In particular, if $P$ is not $\pi_1$--injective then $P$ is
unknotted, $P\cap Q$ is meridional, and $K'$ is primitive in
$H(P)$.
\end{lem}

\begin{proof}
Let $P_1,P_2$ be the closures of the components of $\partial
H(P)\setminus\partial Q_0$. If $Q_0$ is parallel to, say $P_1$,
then $Q$ is isotopic in $X_K$ to $P_1\cup (Q\cap N(P))\subset
N(P)$, which is clearly isotopic to $P$ in $N(P)$ (see
Figure~\ref{birred} and the proof of Lemma~\ref{two}); thus $Q_0$ is
not parallel into $\partial H(P)$.

Let $R$ be the maximal compression body  of $\partial
H(P)=\partial_+ R$ in $H(P)$ (notation as in
\cite{bonahon1,cassongor}). Since $H(P)$ is irreducible and
atoroidal, either $\partial_- R$ is empty and $H(P)$ is a
handlebody or $R$ is a trivial compression body and
$(H(P),\partial H(P))$ is irreducible and atoroidal. As $Q_0$ is
incompressible in $H(P)$, a similar argument shows that either $V$
($W$) is a handlebody or the pair $(V,\partial V)$ ($(W,\partial
W)$, respectively) is irreducible and atoroidal; thus (a) holds.

If $P\cap Q$ is central and either $H(P)$ is a handlebody or $Q_0$
is boundary compressible, then at least one of the circles
$K',l_1,l_2$, where $l_1,l_2$ are the lifts of the simultaneous
centers $P\cap Q$, has a companion annulus by
Lemma~\ref{pants}(a); this contradicts Lemma~\ref{nocomp}, so (b)
now follows from (a).

For part (c), let $\wt{N}(P)=N(P)\cup_{A_K} N(K)=S^3\setminus\intr
H(P)$, the {\it extended} regular neighborhood of $P$ in $S^3$;
notice that $K'\subset\partial\wt{N}(P)$ and, since $P$ has
integral boundary slope, that $A_K$ and $A'_K$ are parallel in
$N(K)$, so $\wt{N}(P)$ and $N(P)$ are homeomorphic in a very
simple way.

Suppose the pair $(V,\partial V)$ is irreducible; without loss of
generality, we may assume $P_1\subset\partial V$ and
$P_2\subset\partial W$. As none of the circles $\partial P_1$
bounds a disk in $N(P)$, the pair of pants $P_1$ is incompressible
in $N(P)$ and hence in $\wt{N}(P)$; thus, since $(V,\partial V)$
is irreducible, it is not hard to see that, for the manifold
$\wt{M}^3=\wt{N}(P)\cup_{P_1\cup (V\cap N(K))} V$ (see
Figure~\ref{birred}), the pair $(\wt{M}^3,\partial\wt{M}^3)$ is
irreducible. As $S^3=\wt{M}^3\cup W$, $\partial W$ must compress
in $W$, so $W$ is a handlebody by (a).

\begin{figure}
\Figw{1.4in}{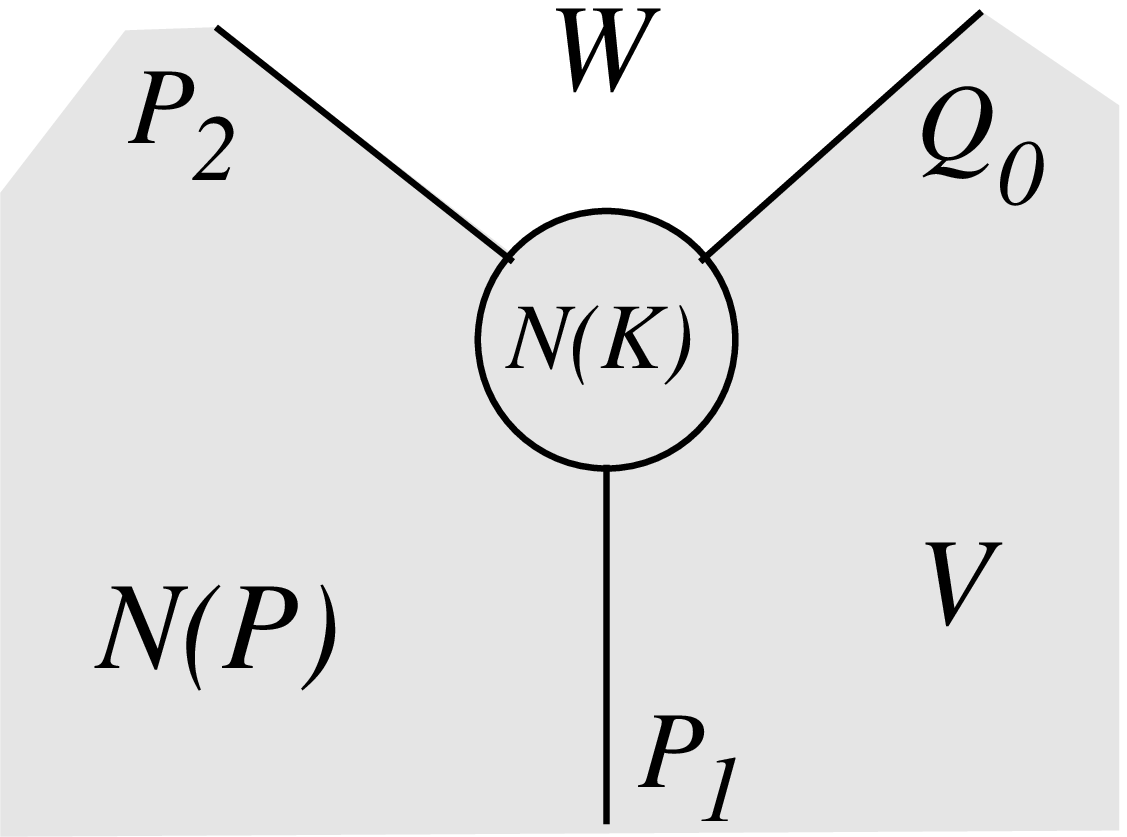}{The manifold
$\wt{M}^3=\wt{N}(P)\cup_{P_1\cup (V\cap N(K))} V$}{birred}
\end{figure}

Finally, if $P$ is not $\pi_1$--injective then $T_P$ compresses in
$H(P)$, so $H(P)$ is a handlebody by (a), $K'$ is primitive in
$H(P)$ by Lemmas~\ref{comp} and \ref{nocomp}, and $P\cap Q$ is
meridional by (b).
\end{proof}

\begin{lem}\label{mer}
If $P,Q,R\in\sk(K,r)$ are distinct elements and each intersection
$P\cap Q,P\cap R$ is central or meridional, then  $P\cap Q$ and
$P\cap R$ are isotopic in $P$.
\end{lem}

\begin{proof}
We will assume that $P\cap Q$ and $P\cap R$ are not isotopic in
$P$ and obtain a contradiction in all possible cases. Since any
two meridian circles of $P$ are isotopic in $P$, we may assume
that $P\cap Q=m$ or $c_1\cup c_2$ and $P\cap R=c'_1\cup c'_2$,
where $m$ is the meridian of $P$ and $c_1,c_2$ and $c'_1,c'_2$ are
two non isotopic pairs of disjoint centers of $P$; we write
$\partial Q_0=\partial Q\cup\alpha_1\cup\alpha_2$ and $\partial
R_0=\partial R\cup l'_1\cup l'_2$, where $\alpha_1,\alpha_2\subset
T_P$ are the lifts $m_1,m_2$ of $m$ or $l_1,l_2$ of $c_1,c_2$, and
$l'_1,l'_2\subset T_P$ are the lifts of $c'_1,c'_2$. Isotope
$P\cap Q$ and $P\cap R$ in $P$ so as to intersect transversely and
minimally; then their lifts $\alpha_1\cup \alpha_2$ and $l_1'\cup
l_2'$ will also intersect minimally in $\partial H(P)$. Finally,
isotope $Q_0,R_0$ in $H(P)$ so as to intersect transversely with
$|Q_0\cap R_0|$ minimal; necessarily, $|Q_0\cap R_0|>0$, and any
circle component of $Q_0\cap R_0$ is nontrivial in $Q_0$ and
$R_0$.

Let $G_{Q_0}=Q_0\cap R_0\subset Q_0$, $G_{R_0}=Q_0\cap R_0\subset
R_0$ be the graphs of intersection of $Q_0,R_0$. Following
\cite{cgls}, we think of the components of $\partial Q_0$ as {\it
fat vertices} of $G_{Q_0}$, and label each endpoint of an arc of
$G_{Q_0}$ with $1'$ or $2'$ depending on whether such endpoint
arises from an intersection involving $l'_1$ or $l'_2$,
respectively; the graph $G_{R_0}$ is labelled with $1,2$ in a
similar way. Such a graph is {\it essential} if each of its
components is essential in the corresponding surface. As $P\cap R$
is central, $R_0$ is boundary incompressible in $H(P)$ by
Lemma~\ref{hbdy-irred}(b) and so $G_{Q_0}$ is essential;
similarly, $G_{R_0}$ is essential if $P\cap Q$ is central. Thus,
if $G_{R_0}$ has inessential arcs then $P\cap Q$ is meridional
and, by minimality of $|Q_0\cap R_0|$, $Q_0$ boundary compresses
along an essential arc of $G_{Q_0}$. By Lemma~\ref{circles}, since
the meridian $m$  and the centers $c'_1,c'_2$ can be isotoped so
that $|m\cap c_j'|=1$, it follows by minimality of $|\partial
Q_0\cap \partial R_0|$ that $|m_i\cap l'_j|=1$ for $i,j=1,2$.
Hence any inessential arc of $G_{R_0}$ has one endpoint in $m_1$
and the other in $m_2$ and so, by Case 1 of Lemma~\ref{pants}(a),
since $Q_0$ is not parallel into $\partial H(P)$ by
Lemma~\ref{hbdy-irred}, $Q_0$ boundary compresses to a companion
annulus of $K'$ in $H(P)$, contradicting Lemma~\ref{nocomp}.

Therefore the graphs $G_{Q_0},G_{R_0}$ are always essential. Since
$P\cap Q$ and $P\cap R$ are not isotopic in $\partial H(P)$, any
circle component of $Q_0\cap R_0$ must be parallel to $\partial
Q,\partial R$ in $Q_0,R_0$, respectively, so, by minimality of
$|Q_0\cap R_0|$, $K'$ must have a companion annulus in $H(P)$,
contradicting Lemma~\ref{nocomp}. Thus $Q_0\cap R_0$ has no circle
components. We consider two cases.

\setcounter{case}{0}
\begin{case}
$P\cap Q=m$
\end{case}

In this case, we have seen that $|m_i\cap l'_j|=1$ for $i,j=1,2$.
Therefore the $m_1,m_2$ and $l'_1,l'_2$ fat vertices of $G_{Q_0}$
and $G_{R_0}$, respectively, each have valence 2, and so, by
essentiality, each graph $G_{Q_0},G_{R_0}$ consist of two parallel
arcs, one annulus face, and one disk face $D$. The situation is
represented in Figure~\ref{graph-r}(a), where only $G_{Q_0}$ is
shown; notice that the labels $1',2'$ must alternate around
$\partial D$.

\begin{figure}
\Figw{4in}{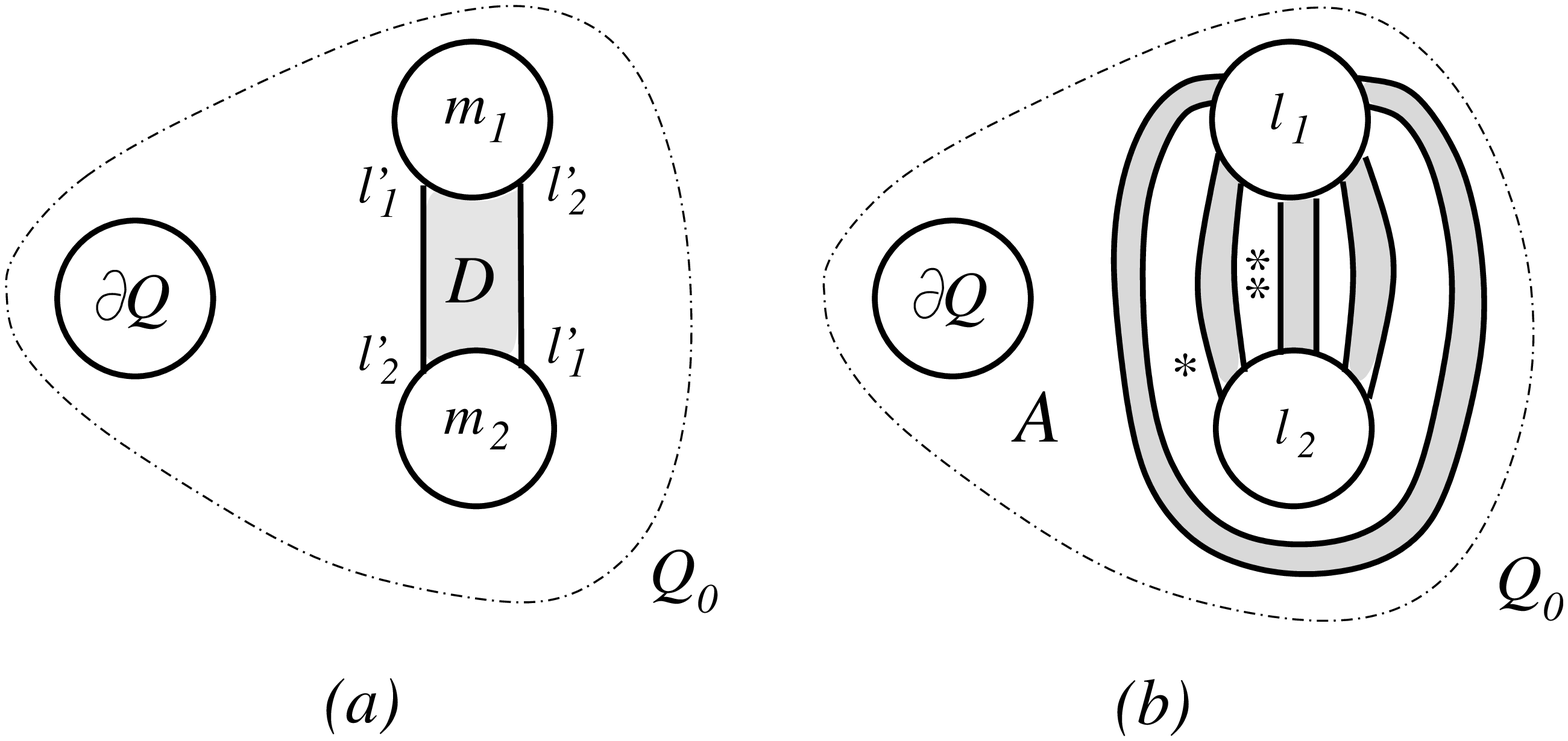}{The possible graphs $G_{Q_0}$}{graph-r}
\end{figure}

Let $V,W$ be the closures of the components of $H(P)\setminus
R_0$, and suppose $Q_0\cap V$ contains the disk face $D$ of
$G_{Q_0}$. By minimality of $|Q_0\cap R_0|$ and $|\partial
Q_0\cap\partial R_0|$, $\partial D$ is nontrivial in $\partial V$,
hence $V$ is a handlebody by Lemma~\ref{hbdy-irred} and $K'$ is
primitive in $V$ by Lemmas~\ref{comp} and \ref{nocomp}. Since the
labels $1',2'$ alternate around $\partial D$, $D$ must be
nonseparating in $V$; thus there is an essential disk $D'$ in $V$
disjoint from $D$ such that $D,D'$ form a complete disk system for
$V$. Now, for $i=1,2$, $|\partial D\cap l'_i|=2$ and so $D\cdot
l'_i$ is even. If the intersection number $D'\cdot l'_i$ is even
for some $i=1,2$ then $l'_i$ is homologically trivial mod 2 in
$V$, while if $D'\cdot l'_i$ is odd for $i=1,2$ then $l'_1\cup
l'_2$ is homologically trivial mod 2 in $V$. Thus one of
$l'_1,l'_2,$ or $l'_1\cup l_2'$ bounds a surface in $V$, hence in
$H(P)$, contradicting Lemma~\ref{lifts}. Therefore this case does
not arise.

\begin{case}
$P\cap Q=c_1\cup c_2$
\end{case}

Recall that $c_1\cup c_2$ and $c_1'\cup c_2'$ intersect
transversely and minimally in $P$, so their lifts $l_1\cup l_2$
and $l_1'\cup l_2'$ also intersect minimally in $\partial H(P)$.
Moreover, after exchanging the roles of $Q_0,R_0$ or relabelling
the pairs $c_1,c_2$ and $c_1',c_2'$, if necessary, we must have
that $|c_1\cap(c'_1\cup c'_2)|=2n+1$ and $|c_2\cap(c'_1\cup
c'_2)|=|2n-1|$ for some integer $n\geq 0$: this can be easily seen
by viewing $P$ as the union of an annulus $\mc{A}$ with core $m$
and a rectangle $\mc{B}$ with core $b$, as in Lemma~\ref{circles},
and isotoping $c_1\cup c_2$ and $c_1'\cup c_2'$ so as to intersect
$m\cup b$ minimally. Figure~\ref{graph-c} represents  two pairs
$c_1\cup c_2$ and $c_1'\cup c_2'$ with minimal intersections; in
fact, all examples of such pairs can be obtained from
Figure~\ref{graph-c}(a) by choosing the rectangle $\mc{B}$
appropriately so that the arcs of $(c_1\cup c_2)\cap \mc{A}$ are
as shown, and suitably Dehn--twisting the pair $c_1'\cup c_2'$
along $m$.

\begin{figure}
\Figw{3.8in}{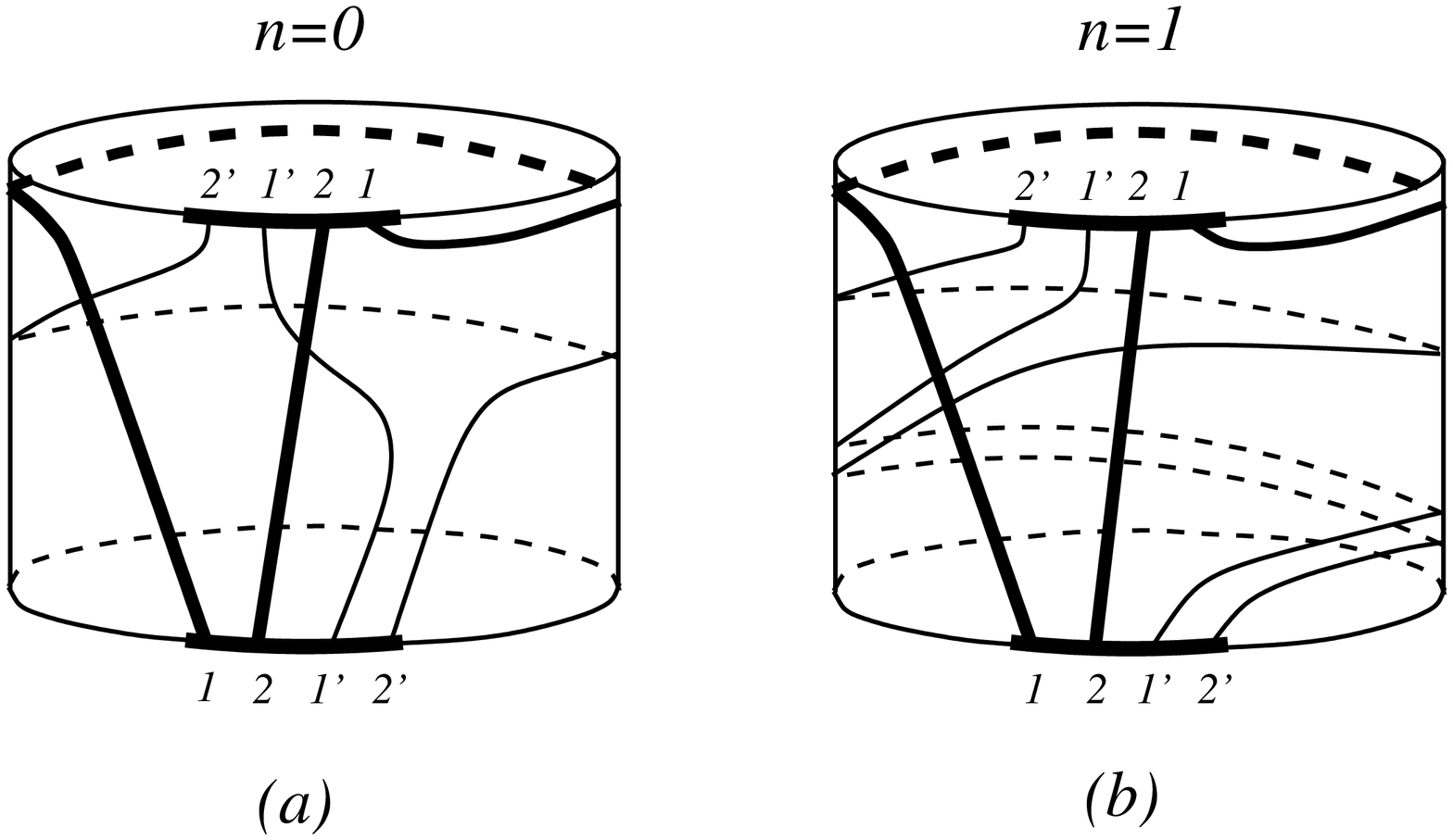}{Minimally intersecting pairs
$c_1\cup c_2,c_1'\cup c_2'$}{graph-c}
\end{figure}

Thus $|l_1\cap(l'_1\cup l'_2)|=4n+2$ and $|l_2\cap(l'_1\cup
l'_2)|=|4n-2|$, and so the essential graphs $G_{Q_0},G_{R_0}$ must
both be of the type shown in Figure~\ref{graph-r}(b) (which is in
fact produced by the intersection pattern of
Figure~\ref{graph-c}(b), where $n=1$). Let $V,W$ be the closures of
the components of $H(P)\setminus R_0$. If $n>1$ then the graph
$G_{Q_0}$ has the two disk components $D_1,D_2$ labelled $*$ and
$**$ in Figure~\ref{graph-r}(b), respectively, as well as an annulus
face $A$ with $\partial_1 A=\partial Q$, all lying in, say, $V$.
Thus $V$ is a handlebody  by Lemma~\ref{hbdy-irred} and, by
minimality of $|\partial Q_0\cap\partial R_0|$ and the
essentiality of $G_{Q_0},G_{R_0}$, the four disjoint circles
$\partial_1 A,\partial_2 A,\partial D_1,\partial D_2$ are all
essential in $\partial V$ and intersect $\partial R_0$ minimally.
Given that $|\partial_1 A\cap (l_1'\cup l_2')|=0$, $|\partial_2
A\cap (l_1'\cup l_2')|=2$, $|\partial D_1\cap (l_1'\cup l_2')|=6$,
and $|\partial D_2\cap (l_1'\cup l_2')|=4$, no two of such four
circles can be isotopic in $\partial V$, an impossibility since
$\partial V$ has genus two. The case $n=0$ is similar to Case 1
and yields the same contradiction.

Finally, for $n=1$ the intersection pattern in $P$ between
$c_1\cup c_2$ and $c_1'\cup c_2'$ must be the one shown in
Figure~\ref{graph-c}(b), and it is not hard to see that only two
labelled graphs $G_{Q_0}$ are produced, up to combinatorial
isomorphism. The first possible labelled graph is shown in
Figure~\ref{graph-ll}(a); capping off $l_1,l_2$ and $l_1',l_2'$ with
the corresponding Moebius bands $Q\cap N(P)$, we can see that
$Q\cap R$ consists of a single circle component (shown in
Figure~\ref{graph-ll}(a) as the union of the broken and solid
lines), which must be a meridian by Lemma~\ref{mlc}, contradicting
Case 1 since $Q\cap P$ is central. In the case of
Figure~\ref{graph-ll}(b), the disk faces $D_1,D_2$ of $G_{Q_0}$ lie
both in, say, $V$, so $V$ is a handlebody by
Lemma~\ref{hbdy-irred} and $K'$ is primitive in $V$ by
Lemmas~\ref{comp} and \ref{nocomp}. Moreover, since $D_1$ and
$D_2$ do not intersect $l'_1\cup l'_2$ in the same pattern, $D_1$
and $D_2$ are nonisotopic compression disks of $\partial
V\setminus K'$ in $V$, so at least one of the disks $D_1,D_2$ must
be nonseparating in $V$. As either disk $D_1,D_2$ intersects each
circle $l'_1,l'_1$ an even number of times, we get a contradiction
as in Case 1. Therefore this case does not arise either.
\end{proof}

\begin{figure}
\Figw{3.8in}{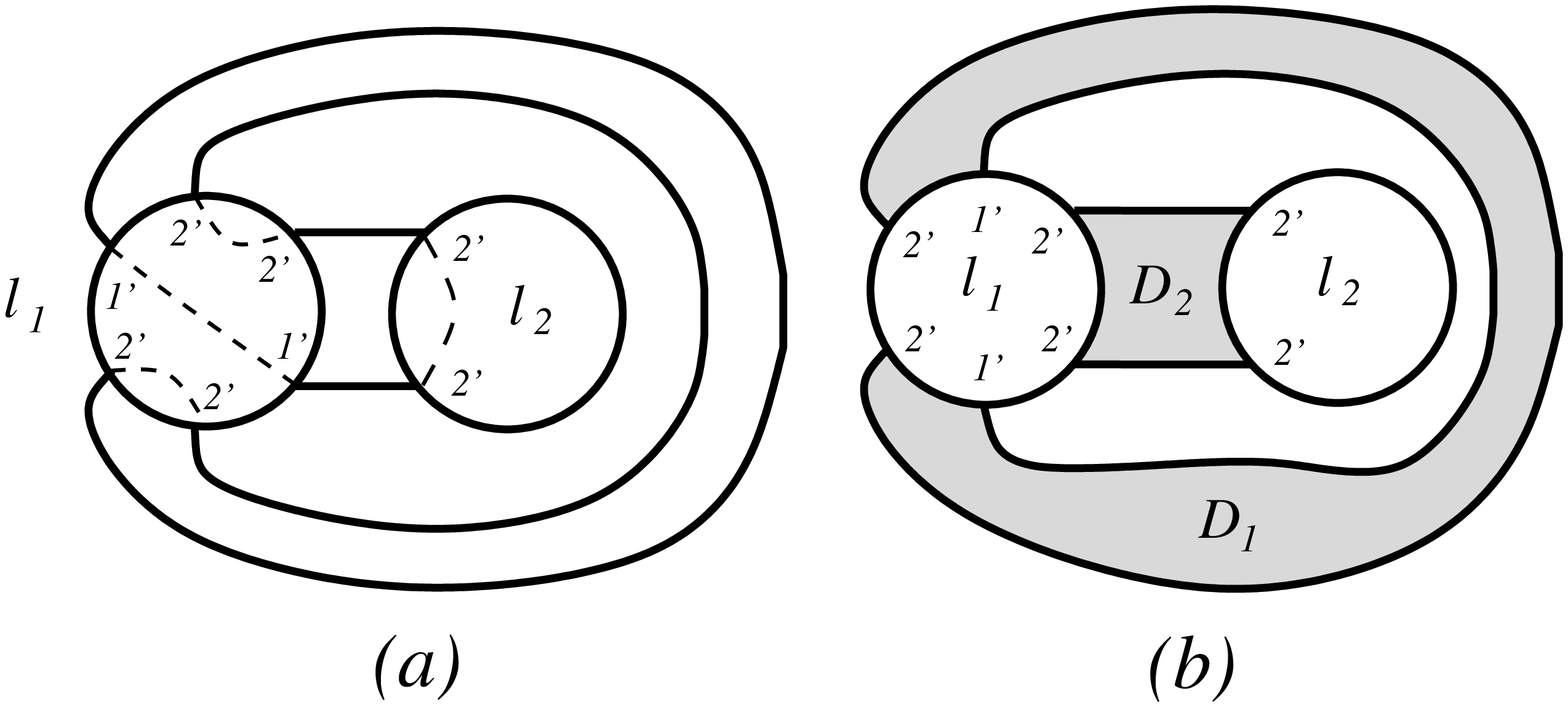}{$G_{Q_0}$ graphs from central
intersections}{graph-ll}
\end{figure}

The next corollary summarizes some of our results so far.

\begin{cor}\label{mer-cen}
If $K$ is hyperbolic and $|\sk(K,r)|\geq 2$ then $\sk(K,r)$ is
either meridional or central; in the latter case, the link
$c_1\cup c_2$ obtained as the intersection of any two distinct
elements of $\sk(K,r)$ is unique in $X_K$ up to isotopy.
\hfill\qed
\end{cor}

If $|\sk(K,r)|\geq 2$ and $P\in\sk(K,r)$, let
$\alpha_1,\alpha_2\subset H(P)$ be the lifts of the common
meridian or pair of disjoint centers of $P$ which, by
Corollary~\ref{mer-cen}, are determined by the elements of
$\sk(K,r)$, and define $P(K,r)$ as the collection of all pairs of
pants $X\subset H(P)$ with $\partial X=K'\cup\alpha_1\cup\alpha_2$
and not parallel into $\partial H(P)$, modulo isotopy. Notice that
$\sk(K,r)\setminus\set{P}$ embeds in $P(K,r)$ by
Corollary~\ref{mer-cen}, and so $|\sk(K,r)|\leq |P(K,r)|+1$. Our
strategy for bounding $|\sk(K,r)|$ in the next lemmas will be to
bound $|P(K,r)|$.

\begin{lem}\label{three}
Let $K$ be a hyperbolic knot with $|\sk(K,r)|\geq 2$. If
$\sk(K,r)$ is central then $|\sk(K,r)|\leq 3$.
\end{lem}

\begin{proof}
Fix $P\in\sk(K,r)$, and suppose $Q_0,R_0,S_0$ are distinct
elements of $P(K,r)$, each with boundary isotopic to $K'\cup
l_1\cup l_2$, where $l_1,l_2$ are the lifts of some fixed pair of
disjoint centers of $P$. Since neither $K',l_1$, nor $l_2$ have
any companion annuli in $H(P)$ by Lemma~\ref{nocomp},
$Q_0,R_0,S_0$ can be isotoped in $H(P)$ so as to become mutually
disjoint. Each of the surfaces $Q_0,R_0,S_0$ separates $H(P)$, and
we may assume that $R_0$ separates $Q_0$ from $S_0$ in $H(P)$.

Let $V,W$ be the closures of the components of $H(P)\setminus
R_0$, with $Q_0\subset V, S_0\subset W$; then $V$, say, is a
handlebody by Lemma~\ref{hbdy-irred}. As  $Q_0$ is not parallel
into $\partial V$ by Lemma~\ref{hbdy-irred}, it follows from
Lemma~\ref{pants}(a) that one of the circles $\partial Q_0$ has a
companion annulus in $V$, hence in $H(P)$, contradicting
Lemma~\ref{nocomp}. Therefore $|P(K,r)|\leq 2$, and so
$|\sk(K,r)|\leq 3$.
\end{proof}

\begin{lem}\label{two}
Let $K$ be a hyperbolic knot with $|\sk(K,r)|\geq 2$. If
$\sk(K,r)$ has an unknotted element $P$ then $\sk(K,r)$ is
meridional, $|\sk(K,r)|=2$, and some lift of the meridian of $P$
has a companion annulus in $H(P)$.
\end{lem}

\begin{proof}
Let $P,Q\in\sk(K,r)$ be distinct elements with $P$ unknotted, and
let $Q_0=Q\cap H(P)$. Then $P\cap Q$ is meridional by
Lemma~\ref{hbdy-irred}(b), so $\sk(K,r)$ is meridional, and $Q_0$
boundary compresses into a companion annulus $A$ of exactly one of
the lifts $m_1$ or $m_2$ of the meridian circle $m$ of $P$ by
Lemmas~\ref{pants}(a) and \ref{nocomp}. Now, if $P_1,P_2$ are the
closures of the components of $\partial H(P)\setminus\partial Q_0$
then, by Lemma~\ref{pants}(b)(ii), $Q_0=P_i\oplus A$ for some
$i=1,2$, so if $A(m)$ is an annulus in $N(P)$ cobounded by
$m_1,m_2$ then $Q$ is equivalent to one of the Seifert Klein
bottles $R_i=P_i\oplus A\cup A(m)$, $i=1,2$. But $R_1$ and $R_2$
are isotopic in $X_K$: this isotopy is described in
Figure~\ref{isotopy}, where a regular neighborhood $N(P)$ of $P$ is
shown (as a box) along with the lifts of the meridian $m$ of $P$
(as two of the solid dots in the boundary of $N(P)$); regarding
$P_1$ as the closure of a component of $\partial
N(P)\setminus(N(A(m)\cup P)$, the idea is to construct $P_1\oplus
A\cup A(m)$, start pushing $P_1$ onto $P$ using the product
structure of $N(P)\setminus A(m)$, and continue until reaching
$P_2$ on the `other side' of $P$; as the annulus $A\cup A(m)$ is
`carried along' in the process, the end surface of the isotopy is
$P_2\oplus A\cup A(m)$. The lemma follows.
\end{proof}

\begin{figure}
\Figw{3.8in}{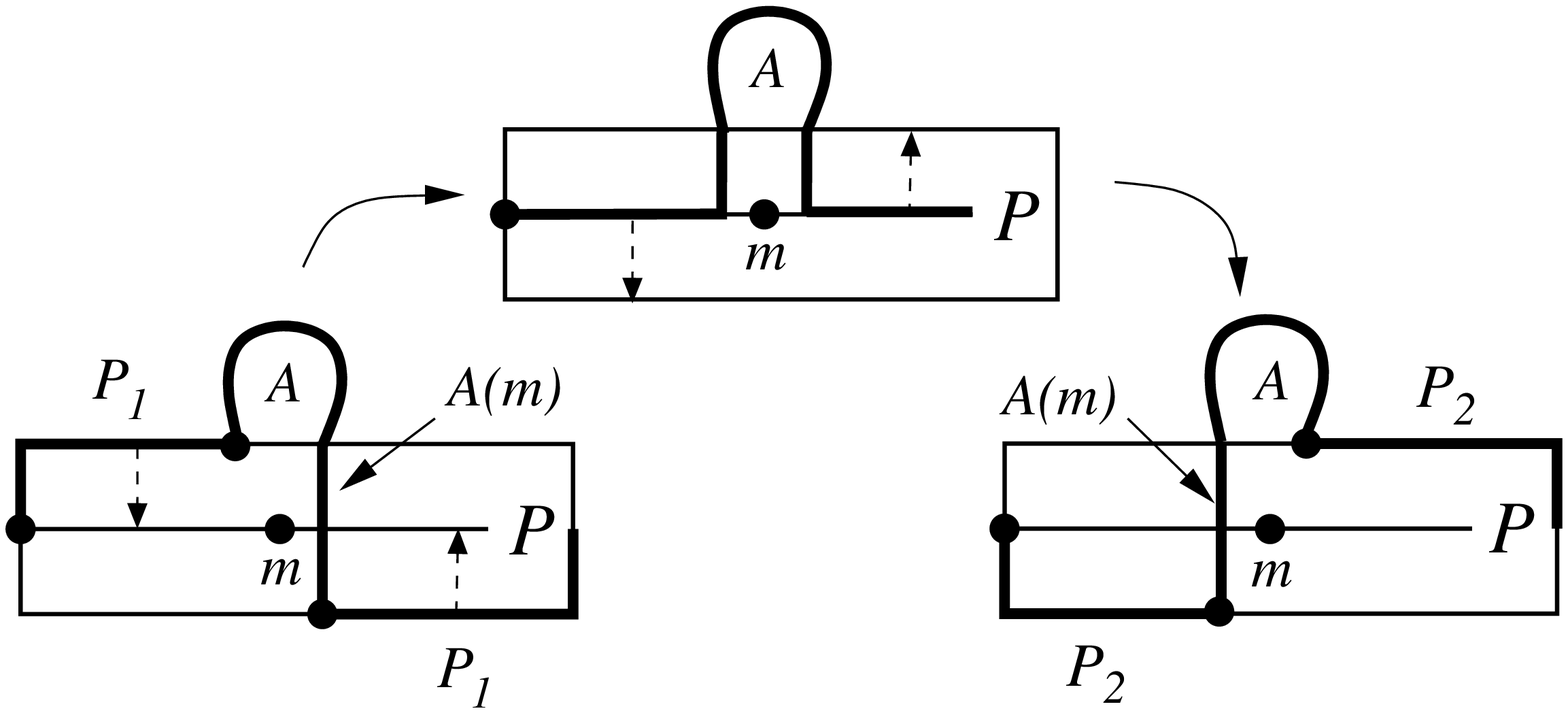}
{An isotopy between $(P_1\oplus A)\cup A(m)$ and $(P_2\oplus
A)\cup A(m)$}{isotopy}
\end{figure}

\begin{rem}
For a hyperbolic knot $K\subset S^3$ with Seifert Klein bottle $P$
and meridian lifts $m_1,m_2$, if there is a pair of pants
$Q_0\subset H(P)$ with $\partial Q_0=K'\cup m_1\cup m_2$ which is
not parallel into $\partial H(P)$, it may still be the case that
the Seifert Klein bottle $Q=Q_0\cup A(m)$ is equivalent to $P$ in
$X_K$. An example of this situation is provided by the hyperbolic
2--bridge knots with crosscap number two; see
Section~\ref{pretzels} and the proof of Theorem~\ref{t2} for more
details.
\end{rem}

\begin{lem}\label{four}
Let $K$ be a hyperbolic knot with $|\sk(K,r)|\geq 2$. If
$\sk(K,r)$ is meridional then $|\sk(K,r)|\leq 6$.
\end{lem}

\begin{proof}
Fix $P\in\sk(K,r)$, and suppose $Q_0,R_0,S_0,T_0$ are distinct
elements of $P(K,r)$ which can be isotoped in $H(P)$ so as to be
mutually disjoint; we may assume $R_0$ separates $Q_0$ from
$S_0\cup T_0$, while $S_0$ separates $Q_0\cup R_0$ from $T_0$. Let
$V,W$ be the closure of the components of $H(P)\setminus R_0$,
with $Q_0\subset V$ and $S_0\cup T_0\subset W$. Then $W$ can not
be a handlebody by Lemma~\ref{pants}(b)(iii), hence, by
Lemmas~\ref{pants}(a) and \ref{hbdy-irred}(c), $V$ is a handlebody
and $Q_0$ gives rise to a companion annulus in $V$ of some lift of
the meridian of $P$. However, a similar argument shows that $T_0$
gives rise to a companion annulus of some lift of the meridian of
$P$ in $W$, contradicting Lemmas~\ref{soltor} and \ref{nocomp}.

Therefore at most three distinct elements of $P(K,r)$ can be
isotoped at a time so as to be mutually disjoint in $H(P)$.
Consider the case when three such disjoint elements exist, say
$Q_0,R_0,S_0$. Let $E,F,G,H$ be the closure of the components of
$H(P)\setminus(Q_0\cup R_0\cup S_0)$, and let $P_1,P_2$ be the
closure of the pair of pants components of $\partial
H(P)\setminus\partial (Q_0\cup R_0\cup S_0)$, as shown
(abstractly) in Figure~\ref{mer1}.

Suppose $T_0\in P(K,r)\setminus\{Q_0,R_0,S_0\}$ has been isotoped
in $H(P)$ so as to intersect $Q_0\cup R_0\cup S_0$ transversely
and minimally. Then $|T_0\cap(Q_0\cup R_0\cup S_0)|>0$ by the
above argument. Since $T_0\cap(Q_0\cup R_0\cup S_0)$ consists of
circles parallel to $\partial T_0$, at least one component of
$\partial T_0$ must have a companion annulus $A_1$ in $H(P)$; the
annulus $A_1$ may be constructed from an outermost annulus
component of $T_0\setminus(Q_0\cup R_0\cup S_0)$, and hence can be
assumed to lie in one of the regions $E,F,G$ or $H$. By
Lemma~\ref{nocomp}, $A_1$ is a companion of some lift $m_1$ of the
meridian of $P$, and by Lemma~\ref{lifts}(b)(c) the components of
$T_0\cap(Q_0\cup R_0\cup S_0)$ are all mutually parallel in each
of $T_0,Q_0,R_0,S_0$. So, if, say, $|T_0\cap Q_0|\geq 2$, then
$m_1$ gets at least one companion annulus on either side of $Q_0$,
which is not possible by Lemma~\ref{soltor}; thus $|T_0\cap X|\leq
1$ for $X=Q_0,R_0,S_0$.

Let $T'_0$ be the closure of the pants component of
$T_0\setminus(Q_0\cup R_0\cup S_0)$; then $T'_0$ lies in one of
$E,F,G$, or $H$, and so $T_0'$  must be isotopic to
$P_1,P_2,Q_0,R_0$, or $S_0$. Also, since $A_1$ is unique up to
isotopy by Lemma~\ref{soltor}, we can write $T_0=T'_0\oplus A_1$.
Therefore, the only choices for $T_0$ are $X\oplus A_1$ for
$X=P_1,Q_0,R_0,S_0$; here we exclude $X=P_2$ since $P_1\oplus A_1$
and $P_2\oplus A_1$ give rise to isotopic once-punctured Klein
bottles in $X_K$ by the proof of Lemma~\ref{two}.

If $A_1\subset E$ (the case $A_1\subset H$ is similar) then
$P_1\oplus A_1=Q_0$ and $Q_0\oplus A_1=P_1$ (see
Figure~\ref{mer1}(a)), hence $T_0=R_0\oplus A_1$ or $S_0\oplus A_1$.
If $A_1\subset F$ (the case $A_1\subset G$ is similar) then
$R_0\oplus A_1=Q_0$ and $Q_0\oplus A_1=R_0$ (see
Figure~\ref{mer1}(b)), hence $T_0=P_1\oplus A_1$ or $S_0\oplus A_1$.
In either case we have $|P(K,r)|\leq 5$, and hence $|\sk(K,r)|\leq
6$. Finally, if at most two distinct elements of $P(K,r)$ can be
isotoped so as to be disjoint in $H(P)$, it is not hard to see by
an argument similar to the above one that in fact the smaller
bound $|\sk(K,r)|\leq 4$ holds.
\end{proof}

\begin{figure}
\Figw{4.5in}{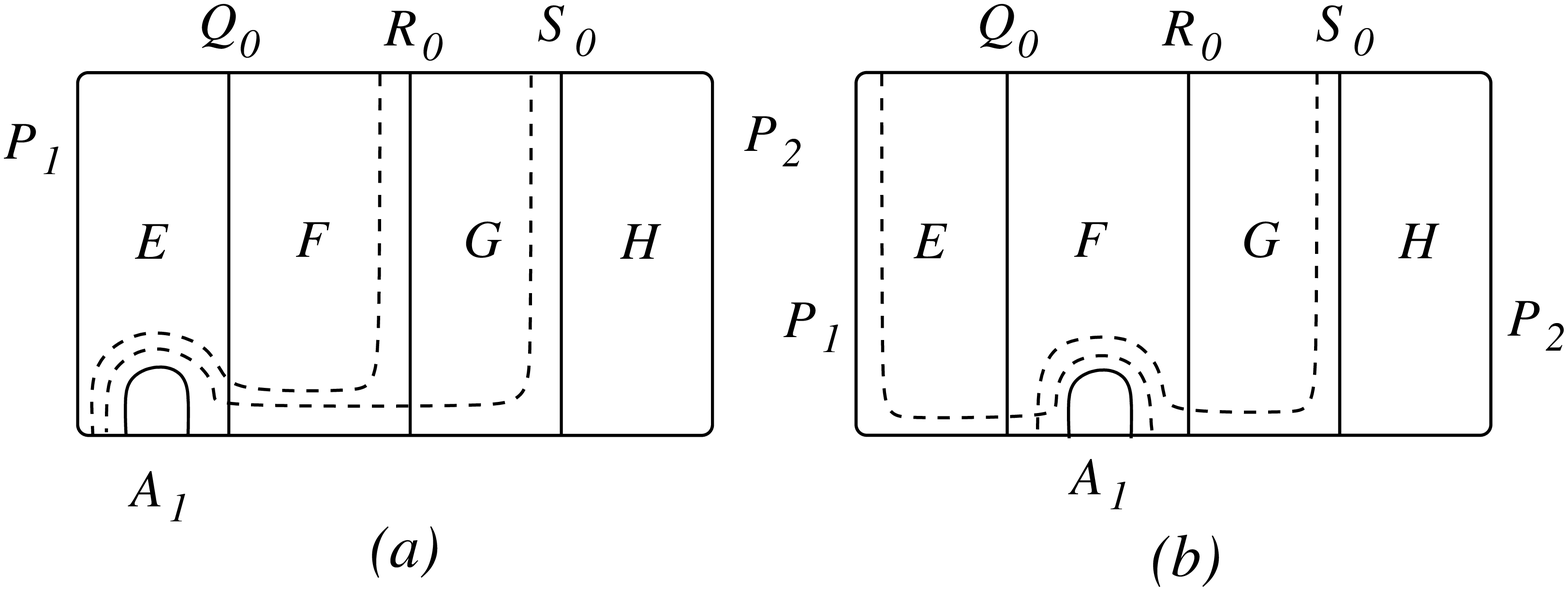}{The possible pairs of pants
$T_0$ (in broken lines) in $H(P)$}{mer1}
\end{figure}

\begin{rem}
It is possible to realize the bound $|P(K,r)|=3$, so
$|\sk(K,r)|\leq 4$, as follows. By \cite[Theorem 1.1]{myers1}, any
unknotted solid torus $S^1\times D^2$ in $S^3$ contains an
excellent properly embedded arc whose exterior $V\subset S^1\times
D^2$ is an excellent manifold with boundary of genus two; in
particular, $(V,\partial V)$ is irreducible, $V$ is atoroidal and
anannular, and $S^3\setminus\intr V$ is a handlebody. Let
$H,Q_0,P$ be the genus two handlebody, pair of pants, and
once-punctured Klein bottle constructed in the remark just after
Lemma~\ref{pants} (see Figure~\ref{example52}). Let $H'$ be the
manifold obtained by gluing a solid torus $U$ to $H$ along an
annulus in $\partial U$ which runs at least twice along $U$ and
which is a regular neighborhood of one of the components of
$\partial Q_0$ which is a lift of the meridian of $P$; since such
a component is primitive in $H$, $H'$ is a handlebody. Finally,
glue $V$ and $H'$ together along their boundaries so that $V\cup
H'=S^3$. Using our results so far in this section it can be proved
(cf proof of Lemma~\ref{mer-hyper}) that $K=\partial P$ becomes
a hyperbolic knot in $S^3$, $H(P)$ contains two disjoint
nonparallel pair of pants with boundary isotopic to $K'$ and the
lifts of a meridian $m$ of $P$, and one of the lifts of $m$ has a
companion annulus in $H(P)$, so $|P(K,r)|=3$. The bound
$|P(K,r)|=5$ could then be realized if $V$ contained a  pair of
pants not parallel into $\partial V$ with the correct boundary.
\end{rem}

\begin{proof}[Proof of Theorem~\ref{t1}]
That $\sk(K,r)$ is either central or meridional and parts (a),(b)
follow from Corollary~\ref{mer-cen} and Lemmas~\ref{hbdy-irred},
\ref{three}, \ref{two}, and \ref{four}.

Let $P\in\sk(K,r)$. If $P$ is not $\pi_1$--injective then $P$ is
unknotted and $K'$ is primitive in $H(P)$ by
Lemma~\ref{hbdy-irred}. Thus, there is a nonseparating compression
disk $D$ of $T_P$ in $H(P)$; since $K'$ has no companion annuli in
$H(P)$, it follows that $N(P)\cup N(D)$ is homeomorphic to $X_K$,
hence $K$ has tunnel number one. Moreover, if $H(P)(K')$ is the
manifold obtained from $H(P)$ by attaching a 2--handle along $K'$,
then $H(P)(K')$ is a solid torus and so
$K(r)=N(\wh{P})\cup_{\partial} H(P)(K')$ is a Seifert fibered
space over $S^2$ with at most three singular fibers of indices
$2,2,n$. As the only such spaces with infinite fundamental group
are $S^1\times S^2$ and $RP^3\# RP^3$, that $\pi_1((K(r))$ is
finite follows from Property R \cite{gabai} and the fact that
$K(r)$ has cyclic integral first homology. Thus (c) holds.

If $P$ is unknotted and $\pi_1$--injective then $T_P$ is
incompressible in $H(P)$, hence, by the 2--handle addition theorem
\cite{cassongor}, the pair $(H(P)(K'),\partial H(P)(K'))$
is irreducible. As $K(r)=N(\wh{P})\cup_{\partial} H(P)(K')$, (d)
follows.
\end{proof}

We discuss now two constructions of crosscap number two hyperbolic
knots. The first construction produces examples of meridional
families $\sk(K,r)$ with $|\sk(K,r)|=2$. The second one gives
examples of knots $K$ and surfaces $P,Q$ in $\sk(K,r)$ which
intersect centrally and such that $|\sk(K,r)|\leq 2$.

\subsection{Meridional families}\label{submer}

It is not hard to produce examples of hyperbolic knots $K$
bounding {\it nonequivalent} Seifert Klein bottles $P,Q$ which
intersect meridionally: for in this case one of the surfaces can
be unknotted and the other knotted, making the surfaces clearly
non isotopic. This is the strategy followed by Lyon in
\cite{lyon1} (thanks to V. N\'u\~nez for pointing out this fact)
to construct nonequivalent Seifert {\it tori} for knots, and his
construction can be easily modified to provide infinitely many
examples of hyperbolic knots $K$ with $|\sk(K,r)|=2$, bounding an
unknotted Seifert Klein bottle and a strongly knotted one along
the same slope.

The construction of these knots goes as follows. As in
\cite{lyon1}, let $V$ be a solid torus standardly embedded in
$S^3$, let $A$ be an annulus embedded in $\partial V$ whose core
is a $(\pm 4,3)$ cable of the core of $V$, and let $A'$ be the
closure of $\partial V\setminus A$. We glue a rectangular band $B$
to $\partial A$ on the outside of $V$, as in Figure~\ref{34torus},
with an odd number of half-twists ($-3$ are shown). Then the knot
$K=\partial (A\cup B)$ bounds the Seifert Klein bottles $P=A\cup
B$ and $Q=A'\cup B$ with common boundary slopes; clearly, $P$ and
$Q$ can be isotoped so that $P\cap Q=A\cap A'$ is a simultaneous
meridian. As in \cite{lyon1}, $P$ is unknotted and $Q$ is knotted;
this is clear since $B$ is a tunnel for the core of $A$ but not
for the core of $A'$. It is not hard to check that if $m_1,m_2$
are the lifts of the meridian of $P$ then $m_1$, say, is a power
(a cube) in $H(P)$ while neither $K',m_2$ is primitive nor a power
in $H(P)$. That $K$ has the desired properties now follows from
the next general result.

\begin{figure}
\Figw{3in}{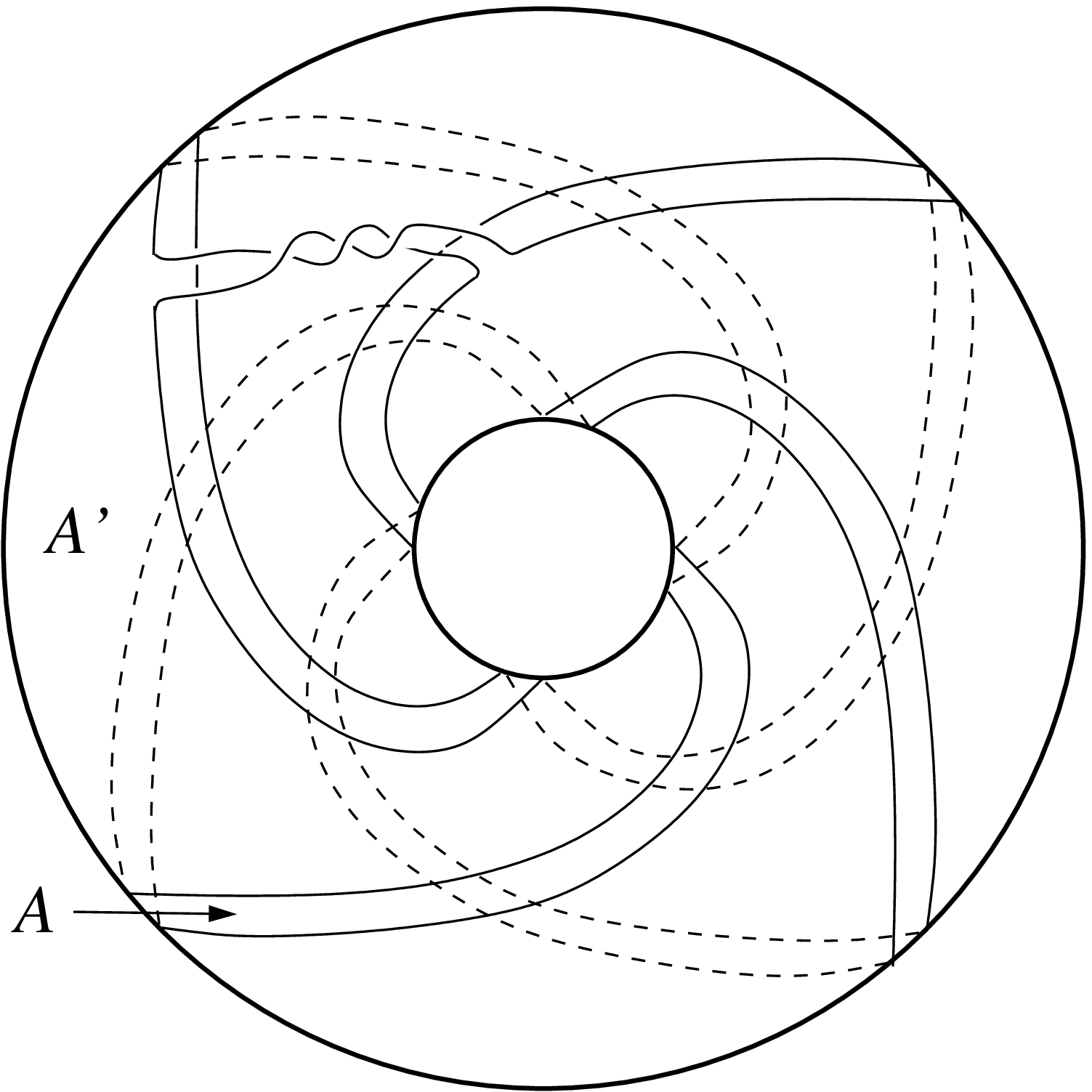}{A modified Lyon's knot}{34torus}
\end{figure}

\begin{lem}\label{mer-hyper}
Let $K$ be a knot in $S^3$ which spans two Seifert Klein bottles
$P,Q$ with common boundary slope $r$, such that $P$ is unknotted,
$Q$ is knotted, and $P\cap Q$ is meridional. For $m_1,m_2$ lifts
of the meridian $m$ of $P$, suppose $m_1$ is a power in $H(P)$ but
neither $K',m_2$ is primitive nor a power in $H(P)$. Then $K$ is
hyperbolic, $P$ is $\pi_1$--injective, and $\sk(K,r)=\set{P,Q}$.
\end{lem}

\begin{proof}
As $K'$ is not primitive nor a power in $H(P)$, $T_P$ is
incompressible in $X_K$ by Lemma~\ref{comp}, so $K$ is a
nontrivial knot and $P$ is $\pi_1$--injective; moreover, $K$ is not
a torus knot by Corollary~\ref{torus} since $Q$ is knotted. Hence
by Lemma~\ref{two} it suffices to show  $X_K$ is atoroidal.

Suppose  $T$ is an essential torus in $X_K$ which intersects $P$
transversely and minimally. Then $P\cap T$ is nonempty and $P\cap
T\subset P$ consists of circles parallel to $\partial P$ and
meridians or longitudes of $P$; by Lemmas~\ref{lifts} and
\ref{nocomp}, since $T$ is not parallel into $\partial X_K$,
it is not hard to see that $P\cap T$ consists of only meridians of
$P$ or only longitudes of $P$. Since the lift of a longitude of
$P$ is also a lift of some center of $P$ then, by
Lemma~\ref{lifts}, either the lift $l$ of some center $c$ of $P$
or both lifts $m_1,m_2$ of the meridian of $P$ have companions in
$H(P)$. The second option can not be the case by Lemma~\ref{comp}
since only $m_1$ is a power in $H(P)$. For the first option,
observe that, since $m$ and $c$ can be isotoped in $P$ so as to
intersect transversely in one point, $l$ can be isotoped in $T_P$
so as to transversely intersect $m_1,m_2$ each in one point.

Suppose $A^*$ is a companion annulus of $l$ in $H(P)$ with
$\partial A^*=l_1\cup l_2$, and let $Q_0=Q\cap H(P)$; notice $Q_0$
is not parallel into $\partial H(P)$, since $Q$ and $P$ are not
equivalent in $X_K$. Isotope $A^*,Q_0$ so as to intersect
transversely and minimally, and let $G_{Q_0}=Q_0\cap A^*\subset
Q_0,G_{A^*}=Q_0\cap A^*\subset A^*$ be their graphs of
intersection; each graph has two arc components. If $G_{Q_0}$ is
inessential then $A^*$ is either parallel into $\partial H(P)$ or
boundary compresses in $H(P)$ into an essential disk disjoint from
$K'$; the first option is not the case, while the latter can not
be the case either by Lemma~\ref{comp} since, by hypothesis, $K'$
is neither primitive nor a power in $H(P)$. If $G_{A^*}$ is
inessential then $Q_0$ boundary compresses in $H(P)$ into a
companion annulus for $K'$, which is also not the case by
Lemma~\ref{comp} since $K'$ is not a power in $H(P)$; for the same
reason, $Q_0\cap A^*$ has no circle components. Thus $G_{Q_0}$ and
$G_{A^*}$ are essential graphs, as shown in Figure~\ref{graph3}.
But, due to the disk face $D_0$ of $G_{Q_0}$, it follows that
$A^*$ runs twice around the solid torus region $R$ cobounded by
$A^*$ and $T_P$; hence $R\cup N(B(c)\subset X_K$ contains a closed
Klein bottle, which is impossible. Thus $X_K$ is atoroidal.
\end{proof}

\begin{figure}\nocolon
\Figw{3.8in}{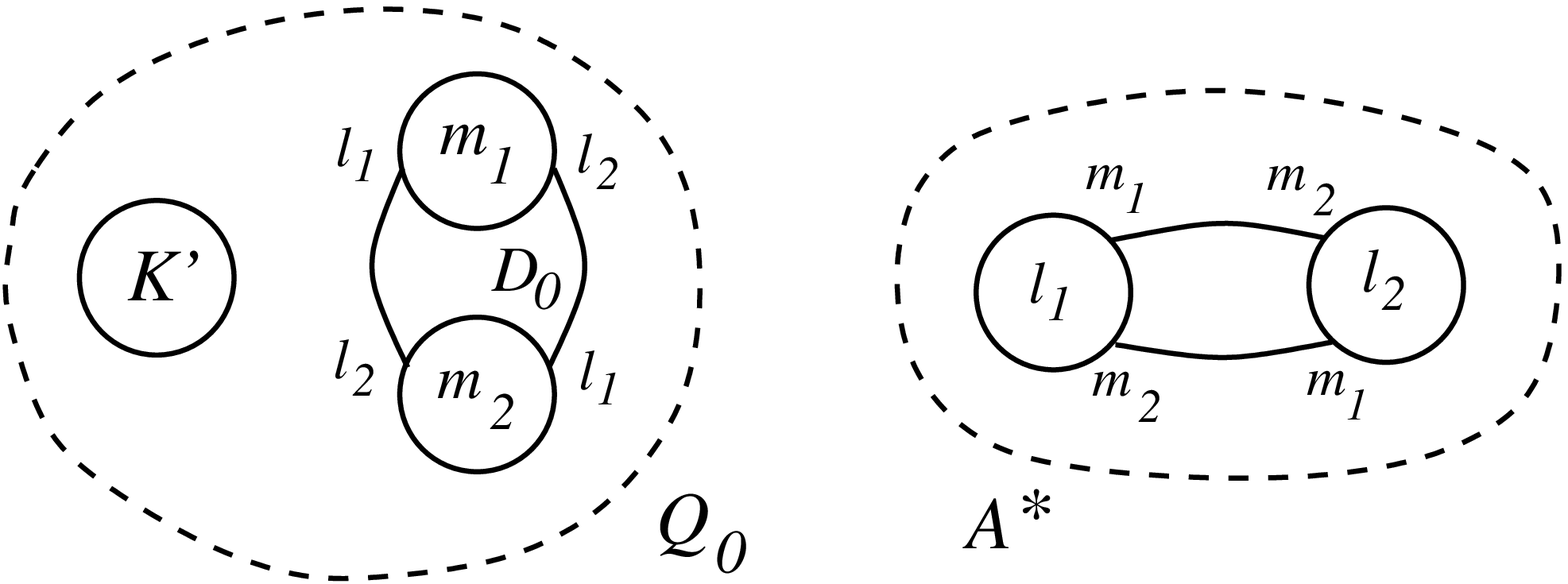}{}{graph3}
\end{figure}

\subsection{Central families}\label{subcen}

Let $S$ be a closed genus two orientable surface embedded in
$S^3$. Suppose there are curves $l_1,l_2\subset S$ which bound
disjoint Moebius bands $B_1,B_2$, respectively, embedded in $S^3$
so that $B_i\cap S=l_i$ for $i=1,2$. Now let $K$ be an embedded
circle in $S\setminus(l_1\cup l_2)$ which is not parallel to
either $l_1$ or $l_2$, does not separate $S$, and separates
$S\setminus(l_1\cup l_2)$ into two pairs of pants, each containing
a copy of both $l_1$ and $l_2$ in its boundary. Let $P_0,Q_0$ be
the closures of the components of $S\setminus(K\cup l_1\cup l_2)$.
Then $K$ bounds the Seifert Klein bottles $P=P_0\cup B_1\cup B_2$
and $Q=Q_0\cup B_1\cup B_2$, which have common boundary slope and
can be isotoped so as to intersect centrally in the cores of
$B_1,B_2$. Any knot $K$ bounding two Seifert Klein bottles $P,Q$
with common boundary slope which intersect centrally can be
constructed in this way, say via the surface $S=(P\asymp Q)\cup
A$, where $A$ is a suitable annulus in $N(K)$ bounded by $\partial
P\cup\partial Q$.

Specific examples can be constructed as follows; however checking
the nonequivalence of two Seifert Klein surfaces will not be as
simple as in the meridional case, as any two such surfaces are
always strongly knotted. Let $S$ be a genus two Heegaard surface
of $S^3$ splitting $S^3$ into genus 2 handlebodies $H,H'$. Let
$l_1,l_2$ be disjoint circles embedded in $S$ which bound disjoint
Moebius bands $B_1,B_2$ in $H$, and let $H_0\subset H$ be the
closure of $H\setminus N(B_1\cup B_2)$. Finally, let $K$ be a
circle in $\partial H\setminus (l_1\cup l_2)$ as specified above,
with $P,Q$ the Seifert Klein bottles induced by $K,l_1,l_2$. It is
not hard to construct examples of $K,l_1,l_2$ satisfying the
following conditions:
\ben
\item[(C1)] $l_1,l_2$ are not powers in $H'$,

\item[(C2)] $K$ is neither primitive nor a power in $H_0,H'$.
\een
The simplest such example is shown in Figure~\ref{central1}; here
$l_1,l_2$ are primitive in $H'$, and $K$ represents
$y^2x^2y^2x^{-2}y^{-2}x^{-2}$, $XYXY^{-1}X^{-1}Y^{-1}$ in
$\pi_1(H),\pi_1(H')$, respectively, relative to the obvious (dual)
bases shown in Figure~\ref{central1}. The properties of $K,P,Q$ are
given in the next result.

\begin{lem}
If $K,l_1,l_2\subset\partial H$ satisfy (C1) and (C2) and $r$ is
the common boundary slope of $P,Q$, then $K$ is hyperbolic, $P$
and $Q$ are strongly knotted, and $|\sk(K,r)|\leq 2$; in
particular, if $P$ and $Q$ are not equivalent then
$\sk(K,r)=\set{P,Q}$.
\end{lem}

\begin{proof}
Observe that $H(P)=H_0\cup_{Q_0} H'$ and $H(Q)=H_0\cup_{P_0} H'$.
That $Q_0$ is incompressible and boundary incompressible in $H(P)$
follows from (C1) and (C2) along with the fact that $l_1,l_2$ are
primitive in $H_0$. Thus $(H(P),\partial H(P))$ is irreducible, so
$K$ is nontrivial and not a cable knot by Theorem~\ref{t4}, and
$P$ (similarly $Q$) is strongly knotted. That $K$ is hyperbolic
follows now from an argument similar to that of the proof of
Lemma~\ref{mer-hyper} and, since both $H_0$ and $H'$ are
handlebodies, the bound $|\sk(K,r)|\leq 2$ follows from the proof
of Lemma~\ref{three}.
\end{proof}

\begin{figure}
\Figw{3.8in}{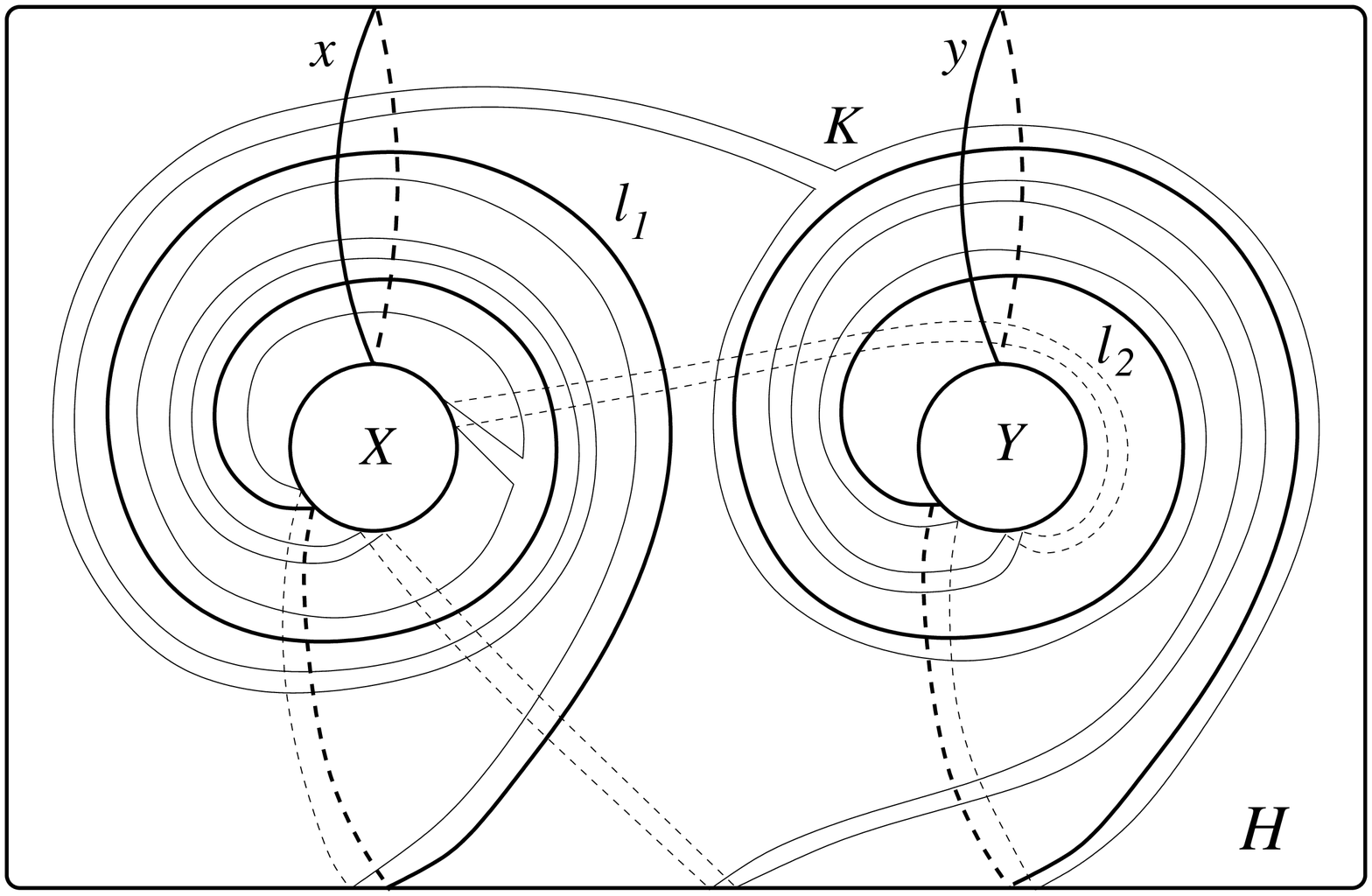}{The circles $K,l_1,l_2$
in $\partial H$}{central1}
\end{figure}

\section{Pretzel knots}\label{pretzels}

We will denote a pretzel knot of length three with the standard
projection shown in Figure~\ref{pretzel2} by $p(a,b,c)$, where the
integers $a,b,c$, exactly one of which is even, count the number
of signed half-twists of each tangle in the boxes. It is not hard
to see that if $\set{a',b',c'}=\set{\varepsilon a,\varepsilon
b,\varepsilon c}$ for $\varepsilon=\pm 1$ then $p(a,b,c)$ and
$p(a',b',c')$ have the same knot type.  For any pretzel knot
$p(a,b,c)$ with $a$ even, the black surface of its standard
projection shown in Figure~\ref{pretzel2} is an algorithmic Seifert
Klein bottle with meridian circle $m$, which has integral boundary
slope $\pm 2(b+c)$ by Lemma~\ref{slope}; an algorithmic Seifert
surface is always unknotted. By \cite{valdez6}, with the exception
of the knots $p(2,1,1)$ (which is the only knot that has two
algorithmic Seifert Klein bottles of distinct slopes produced by
the same projection diagram) and $p(-2,3,7)$, this is the only
slope of $p(a,b,c)$ which bounds a Seifert Klein bottle. Finally,
if at least one of $a,b,c$ is $\pm 1$ then $p(a,b,c)$ has bridge
number at most 2, the only pretzels $p(a,b,c)$ with
$|a|,|b|,|c|\geq 2$ which are torus knots are $p(-2,3,3)$ and
$p(-2,3,5)$, and if one of $a,b,c$ is zero then $p(a,b,c)$ is
either a 2--torus knot or a connected sum of two 2--torus knots.

\begin{figure}
\Figw{1.5in}{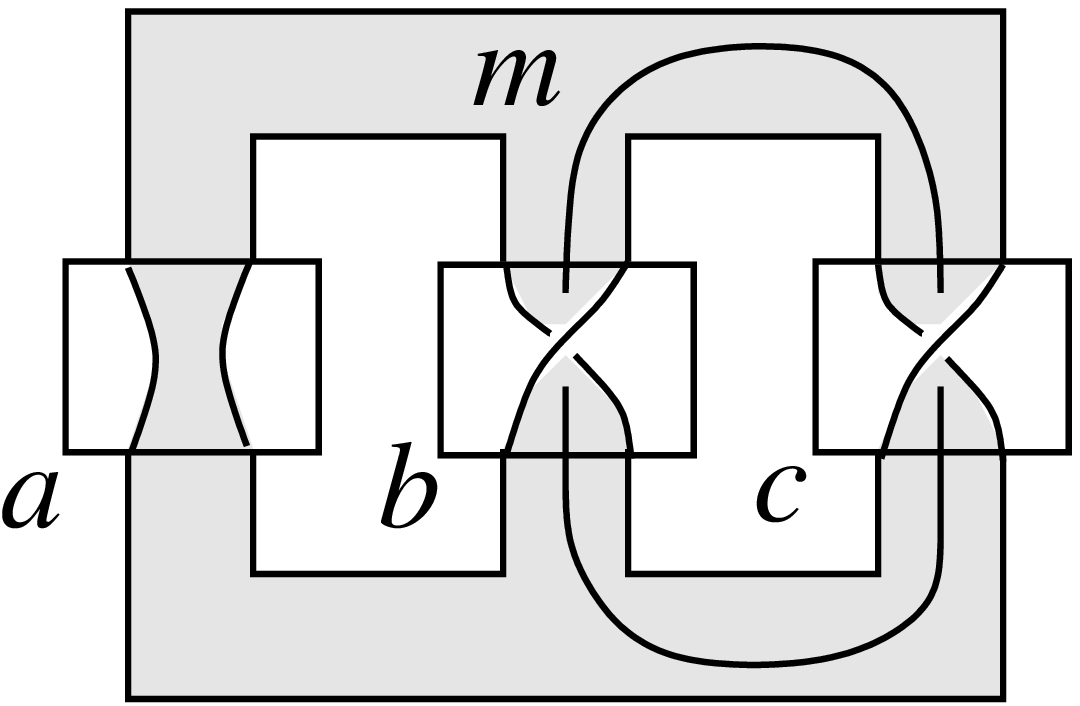}{The pretzel knot $p(a,b,c)$
for $a=0,b=1,c=1$}{pretzel2}
\end{figure}

Now let $F$ be the free group on $x,y$. If $w$ is a cyclically
reduced word in $x,y$ which is primitive in $F$ then, by
\cite{cohen1} (cf \cite{fico1}), the exponents of one of
$x$ or $y$, say $x$, are all $1$ or all $-1$, and the exponents of
$y$ are all of the form $n,n+1$ for some integer $n$. Finally, a
word of the form $x^my^n$ is a proper power in $F$ iff
$\set{m,n}=\set{0,k}$ for some $|k|\geq 2$.

\begin{proof}[Proof of Theorem~\ref{t2}]
Let $K$ be the hyperbolic pretzel knot $p(a,b,c)$ for some
integers $a,b,c$ with $a$ even and $b,c$ odd. Let $P$ be the
unknotted Seifert Klein bottle spanned by $K$ in its standard
projection.

We will show that $K'$ is never primitive in $H(P)$, so $P$ is
$\pi_1$--injective by Lemma~\ref{hbdy-irred}; thus $K(r)$ is
irreducible and toroidal by Theorem~\ref{t1}(d). We will also show
that whenever a lift $m_1,m_2$ of the meridian $m$ of $P$ is a
power in $H(P)$ then $K$ is a 2--bridge knot; in such case, by
\cite[Theorem 1]{thurs4}, $P$ is obtained as a plumbing of an
annulus and a Moebius band (cf \cite{valdez8}) and $P$ is unique
up to isotopy. Along with Lemma~\ref{two}, it will then follow
that $|\sk(K,r)|=1$ in all cases.

The proof is divided into cases, depending on the relative signs
of $a,b,c$. Figure~\ref{pretzel1} shows the extended regular
neighborhood $\wt{N}(P)=N(P)\cup_{A_K} N(K)$ of $P$, which is a
standard unknotted handlebody in $S^3$, along with the circles
$K',m_1,m_2$ with a given orientation. The disks $D_x,D_y$ shown
form a complete disk system for $H(P)$, and give rise to a basis
$x,y$ for $\pi_1(H(P))$, oriented as indicated by the head (for
$x$) and the tail (for $y$) of an arrow. Figure~\ref{pretzel1}
depicts the knot $p(2,3,3)$ and illustrates the general case when
$a,b,c> 0$; we will continue to use the same figure, with suitable
modifications, in all other cases. By the remarks at the beginning
of this section, the following cases suffice.

\setcounter{case}{0}
\begin{case}
$a=2n>0, \ b=2p+1>0, \ c=2q+1>0$
\end{case}

In this case, up to cyclic order, the words for $K',m_1,m_2$ in
$\pi_1(H(P))$ are:
\[\begin{array}{rl}
K'& = y^{q+1}(xy)^px^{n+1}(yx)^py^{q+1}x^{-n}\\
m_1 & =y^q(yx)^{p+1}\\
m_2 & =y^{q+1}(yx)^p
\end{array}
\]
Since $n>0$ and $p,q\geq 0$, the word for $K'$ is cyclically
reduced and both  $x$ and $x^{-1}$  appear in $K'$; thus $K'$ is
not primitive in $H(P)$.

Consider now $m_1$ and $m_2$; as $y$ and $yx$ form a basis of
$\pi_1(H(P))$, if $m_1$ is a power in $H(P)$ then $q=0$, while if
$m_2$ is a power then $p=0$. In either case $K$ is a pretzel knot
of the form $p(\cdot,\cdot,1)$, hence $K$ is a 2--bridge knot.

\begin{figure}\nocolon
\Figw{2.5in}{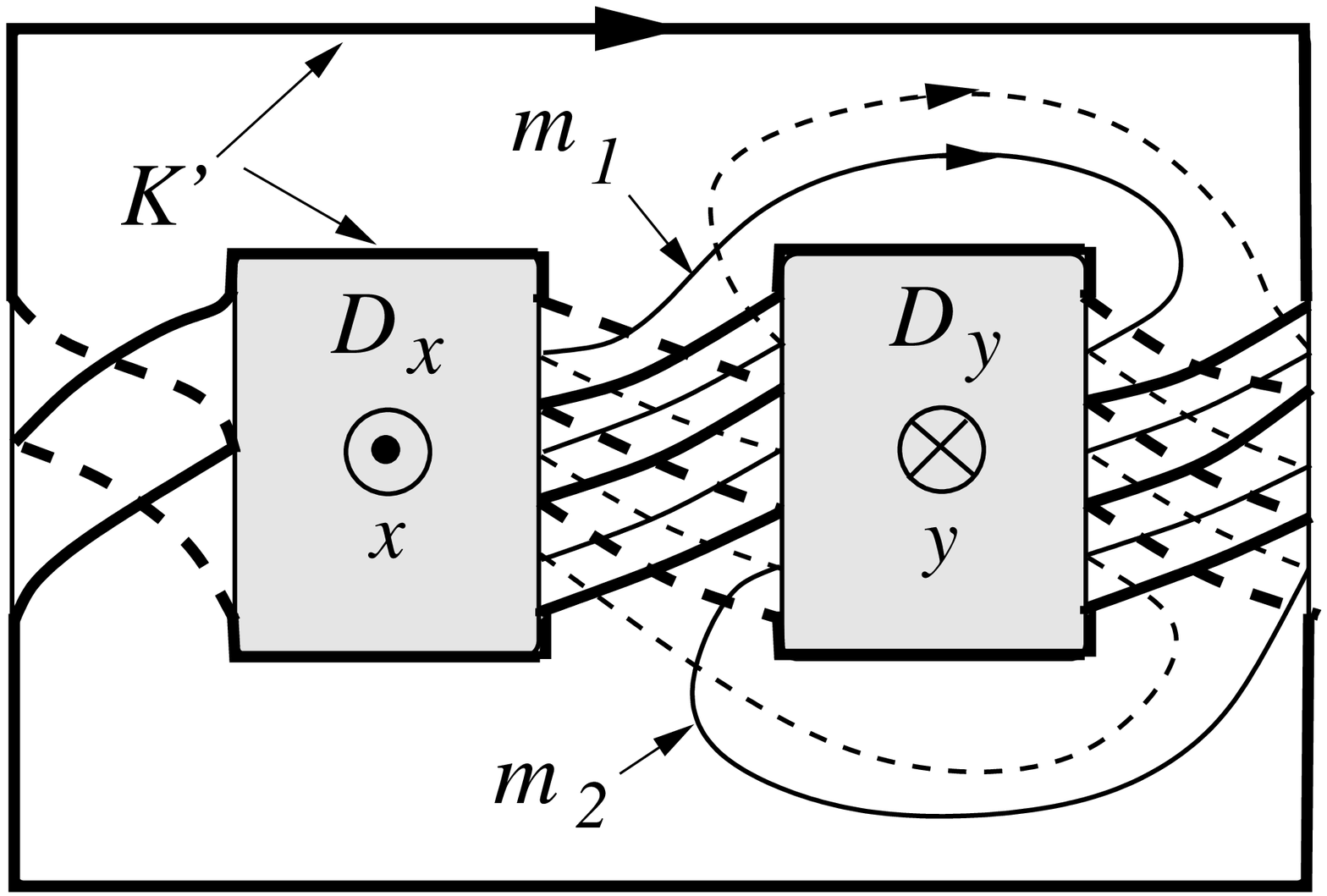}{}{pretzel1}
\end{figure}

\begin{case}
$a=-2n<0, \ b=2p+1>0, \ c=2q+1>0$
\end{case}

This time words for $K',m_1,m_2$ are, up to cyclic order,
\[\begin{array}{rl}
K'& = y^{q+1}(xy)^px^{1-n}(yx)^py^{q+1}x^{n}\\
m_1 & =y^q(yx)^{p+1}\\
m_2 & =y^{q+1}(yx)^p.
\end{array}
\]
If $n>1$ then the word for $K'$ is cyclically reduced and both $x$
and $x^{-1}$ appear in $K'$ and so $K'$ is not primitive in
$H(P)$. If $n=1$ then, switching to the basis $y,u=xy$ of
$\pi_1(H(P))$, $K'$ is represented up to cyclic order by the word
$y^qu^pyu^py^qu$. Observe that if $p=0$ or $q=0$ then
$K=p(-2,1,\cdot)$ which is a 2--torus knot, so $p,q>0$. Thus, if
$K'$ is primitive then necessarily $\set{p,q}=\set{1,2}$ and so
$K=p(-2,3,5)$ is a torus knot. Therefore, $K'$ is not primitive in
$H(P)$. The analysis of the words $m_1$ and $m_2$ is identical to
that of Case 1 and yields the same conclusion.

\begin{case}
$a=2n>0, \ b=2p+1>0, \ c=-(2q+1)<0$
\end{case}

Up to cyclic order, words for $K',m_1,m_2$ are:
\[\begin{array}{rl}
K'& = y^{q}(xy^{-1})^px^{n+1}(y^{-1}x)^py^{q}x^{-n}\\
m_1 & =y^{q+1}(xy^{-1})^{p+1}\\
m_2 & =y^{q}(y^{-1}x)^p.
\end{array}
\]
If $p,q>0$ then the word for $K'$ is cyclically reduced  and
contains all of $x,x^{-1},y,y^{-1}$, so $K'$ is not primitive in
$H(P)$. If $p=0$ then $K'=y^qx^{n+1}y^qx^{-n}$, so $K'$ is
primitive iff $q=0$, in which case $K=p(2n,1,-1)$ is a 2--torus
knot. The case when $q=0$ is similar, therefore $K'$ is not
primitive in $H(P)$.

In this case $m_1$ can not be a power in $H(P)$ for any values of
$p,q\geq 0$, while if $m_2$ is a power then $p=0$ or $q=0$ and
hence $K$ is a 2--bridge knot.
\end{proof}

\Addresses\recd

\end{document}